\documentclass[a4paper]{elsarticle}
\usepackage[utf8x]{inputenc}
\usepackage[english]{babel,isodate}
\usepackage{pgfplotstable}
\usepackage{algorithm,algorithmic}
\usepackage{graphicx}
\usepackage{grffile}
\usepackage{amsmath,amsthm}
\usepackage{amssymb}
\usepackage{todonotes}
\usepackage[integrals]{wasysym}
\usepackage{rotating}
\usepackage{bibgerm}
\usepackage{hyperref}
\usepackage{environ}
\usepackage{enumerate}
\usepackage{float}
\restylefloat{table}
\usepackage[titletoc]{appendix}
\usepackage{tikz}
\usetikzlibrary{patterns}
\usetikzlibrary{decorations.pathreplacing}
\usetikzlibrary{calc,matrix,positioning}
\usepgfplotslibrary{groupplots}
\usepackage{nicefrac}
\usepackage{caption}
\usepackage{fullpage}

\theoremstyle{definition}
\newtheorem{theorem}{Theorem}
\newtheorem{lemma}{Lemma}
\newtheorem*{theorem*}{Theorem}

\newcommand{\dl}{\mathrm{d}}
\newcommand{\interior}{\mathrm{int}}
\newcommand{\exterior}{\mathrm{ext}}

\newcommand\todoblock[2][]{\todo[inline,caption={},#1]{#2}}
\NewEnviron{todoenv}[1][]{\todoblock[#1]{\BODY}}

\newcommand\THcube[2]{$\mathcal{Q}_{#1}/\mathcal{Q}_{#2}$}
\newcommand\RT{\text{RT}}

\newcommand\normal{n}
\newcommand{\tnabla}\nabla

\title{A high-order discontinuous Galerkin pressure robust splitting scheme for incompressible flows
  \\[5mm]
  \small{Dedicated to Mary F.~Wheeler in honour of her 80th birthday}}
\author[heidelberg]{Marian~Piatkowski\corref{cor1}\fnref{fn1}}
\ead{marian.piatkowski@iwr.uni-heidelberg.de}
\author[heidelberg]{Peter~Bastian}
\ead{peter.bastian@iwr.uni-heidelberg.de}
\cortext[cor1]{Corresponding author}
\fntext[fn1]{email: \texttt{marian.piatkowski@iwr.uni-heidelberg.de}}
\address[heidelberg]{Interdisciplinary Center for Scientific Computing, Heidelberg University, Im Neuenheimer Feld 205, 69120 Heidelberg, Germany}

\begin{document}
\begin{abstract}
  The accurate numerical simulation of high Reynolds number incompressible flows is a challenging
  topic in computational fluid dynamics. Classical inf-sup stable methods like the Taylor-Hood
  element or only \(L^2\)-conforming discontinuous Galerkin (DG) methods relax the divergence
  constraint in the variational formulation. However, unlike divergence-free methods, this
  relaxation leads to a pressure-dependent contribution in the velocity error which is proportional
  to the inverse of the viscosity, thus resulting in methods that lack pressure robustness
  and have difficulties in preserving structures at high Reynolds numbers.
  The present paper addresses the discretization of the incompressible Navier-Stokes
  equations with high-order DG methods in the framework of projection methods. The major focus in
  this article is threefold:
  i) We present a novel postprocessing technique in the projection step of the splitting scheme
  that reconstructs the Helmholtz flux in \(H(\text{div})\). In contrast to the previously introduced
  \(H(\text{div})\) postprocessing technique, the resulting velocity field is pointwise
  divergence-free in addition to satisfying the discrete continuity equation.
  ii) Based on this Helmholtz flux \(H(\text{div})\) reconstruction, we obtain a high order in space,
  pressure robust splitting scheme as numerical experiments in this paper demonstrate.
  iii) With this pressure robust splitting scheme, we demonstrate that a robust DG method for
  underresolved turbulent incompressible flows can be realized.
\end{abstract}
\begin{keyword}
  incompressible Navier-Stokes, projection method, high-order discontinuous Galerkin,
  pressure robust method, divergence-free properties, turbulence modeling, implicit LES
\end{keyword}
\maketitle

\section{Introduction}
\label{sec:introduction}

Navier-Stokes fluid flow at high Reynolds numbers plays an increasing role also for the detailled
understanding of subsurface flow processes.
Our original interest was in contributing to the understanding of the surface renewal effect at the interface
of subsurface and atmosphere \cite{KATUL2006117,Piatkowski2019} which is attributed to coherent
turbulent structures interacting with the atmospheric boundary layer at the surface. Other applications include Navier-Stokes
flow in rock fractures \cite{doi:10.1029/162GM07,doi:10.1029/2007GL030545}, flow in wellbores \cite{doi:10.1177/1687814019873250}
as well as flows on the pore-scale \cite{RAEINI20125653,doi:10.1002/fld.3653}.

The numerical simulation of flows at high Reynolds numbers is still a challenging task. Most classical finite
element methods relax the divergence constraint and only enforce the condition
weakly. Recently, it has been understood that the relaxation of the incompressibility constraint
introduces a pressure-dependent contribution in the velocity error for inf-sup stable methods like
the Taylor-Hood element or only \(L^2\)-conforming DG methods. This insight has led to distinguish between
pressure robust and non pressure robust space discretizations, respectively, where a method is called
pressure robust if the velocity error is independent of the continuous pressure. The study of
pressure robust methods is a current, active research topic
\cite{GAUGER201989,LINKE2019350,Schroeder2018,AhmedLinkeMerdon18,doi:10.1137/17M1112017}.
It is connected to the Helmholtz decomposition of vector fields and consequently, a fundamental
invariance property of the incompressible Navier-Stokes equations if the boundary conditions do not
depend on the pressure: Let \((v,p)\) be the solution of the equations with right-hand side \(f\),
then a transformation of the body force \(f\mapsto f+\tnabla\psi\) changes the Navier-Stokes solution
\((v,p)\mapsto(v,p+\psi)\), i.e. the velocity field does not change and the additional forcing is balanced by
the pressure gradient. A desirable property of a discretization scheme is to maintain an unchanged
velocity field under such irrotational force translation. As a simple example let us consider the
stationary Stokes equations
\begin{align*}
  -\nu\Delta v + \tnabla p &= f \hspace*{0.9em} \textup{in} \; \Omega \\
  \nabla\cdot v &= 0 \quad \textup{in} \; \Omega \\
  v &= g \quad \textup{on} \; \partial\Omega \\
  \int_\Omega p \dl x &= 0 \; .
\end{align*}
Employing a \(k\)-th order mixed finite element DG discretization, the following error estimate in
\cite{Girault_adiscontinuous,girault_riviere} was shown:
\begin{equation}
  \label{eq:classical_velocityerror_estimate}
  \bigl\lbrack\! \| v-v_h \| \!\bigr\rbrack
  \leq C h^k \left( |v|_{k+1} + \frac{1}{\nu} |p|_k \right)
\end{equation}
with constant \(C\) independent on \(h\) and \(\nu\), and with the mesh-dependent norm
\(\bigl\lbrack\! \| \cdot \| \!\bigr\rbrack\) for the Stokes problem.
The right-hand side in the a-priori estimate (\ref{eq:classical_velocityerror_estimate}) is a standard
bound for classical mixed finite element methods. Observe that the velocity error depends on the
pressure scaled by the inverse of viscosity. This has the consequence that the velocity
approximation is not robust with respect to irrotational force translations that would solely change
the pressure in the continuous case and therefore this method is not pressure robust. In fact, a
consistency error is hidden in the discrete Helmholtz projector such that discretely divergence-free
functions need not be \(L^2\)-orthogonal to the irrotational fields. The consistency error is
reflected in the pressure term \(\frac{1}{\nu}\left|p\right|_k\) and renders such space
discretizations inaccurate at handling correctly large irrotational parts in a flow, or preserving
flow structures, for high Reynolds numbers, especially in an underresolved turbulence computation.
In contrast, the pressure-dependent term on the right-hand side disappears if a divergence-free
method is employed.

It is thus important for a space discretization to preserve
physical properties at the discrete level. To name only a few, inf-sup stable divergence-free
mixed methods \cite{10.2307/4100078}, 
inf-sup stable \(H(\text{div})\) conforming DG methods \cite{Cockburn2007,Lehrenfeld2010},
inf-sup stable finite element methods with appropriately modified velocity test functions
\cite{LINKE2012837,LINKE2014782,DIPIETRO2016175,doi:10.1137/15M1047696,refId0,LINKE2016304},
and inf-sup stable
\(L^2\)-conforming DG methods with additional consistent stabilization terms
\cite{KRANK2017634,FEHN2018667} have been found lately in order to realize this.

High-order DG methods are increasingly important for a wide range of applications including
computational fluid dynamics given recent development in modern computer architectures. Their block
structure, compact stencils, and high ratio of computation to communication make them well-suited on
modern, many-core, memory-constrained architectures. To this extent, there has been recent work on
high-order DG methods for the incompressible Navier-Stokes equations, see
e.g. \cite{Schroeder2019,doi:10.1002/fld.4763} for coupled solution approaches,
\cite{2016arXiv161200657P}, and \cite{KRANK2017634,FEHN2017392,doi:10.1002/fld.4683} for
discretizations using projection methods.
In \cite{2016arXiv161200657P} we have presented a postprocessing
technique in the Helmholtz projection step based on \(H(\text{div})\) reconstruction of the pressure
correction. The obtained velocity field has been shown to satisfy the discrete continuity
equation so that the procedure defines a discrete projection operator.

Extending this solution approach, we present in this work a novel postprocessing technique that
reconstructs the Helmholtz flux \(w_h - \tnabla_h\psi_h\) in the Raviart-Thomas space where \(w_h\)
denotes the tentative velocity and \(\tnabla_h\psi_h\) is the irrotational correction. The resulting
velocity field satisfies the discrete continuity equation, and is also pointwise divergence-free,
the latter property not being satisfied by the previously introduced \(H(\text{div})\) postprocessing
technique. The reconstruction picks up the notion of a divergence-preserving reconstruction operator
for discontinuous velocity and pressure ansatz spaces. As numerical experiments demonstrate,
exactly the improvement on the pointwise divergence renders the approach to be pressure robust like a
divergence-free mixed method. We thereby substantiate the recently made statement
\cite{GAUGER201989} that the need for pressure robustness emanates from an improved
understanding of mixed methods and the divergence constraint in incompressible flows.
\paragraph{Outline}
The article is organized as follows: In Section \ref{sec:disc-galerk-spat} we introduce the notation
for the discontinuous Galerkin formulation and describe the DG discretization of the incompressible
Navier-Stokes equations. In Section \ref{sec:splitting-method} we introduce
splitting methods where we concentrate on the pressure correction scheme in rotational form (RIPCS)
throughout this work. Then, several variants of discrete Helmholtz decompositions in
the framework of projection methods are reviewed and later compared in the numerical experiments.
In Section \ref{sec:helmh-flux-ravi} we present the novel
postprocessing technique in the projection step of the splitting scheme that reconstructs the
Helmholtz flux \(w_h - \tnabla_h\psi_h\) in \(H(\text{div})\). It is shown that the resulting velocity
field satisfies the discrete continuity equation and is pointwise divergence-free. We assess the
numerical conservation properties and temporal convergence of the obtained splitting scheme by
benchmarks that have been used in a previous publication \cite{2016arXiv161200657P}. In Section
\ref{sec:numerical-results} the aforementioned variants of discrete Helmholtz decomposition are
tested and compared within the pressure correction scheme. The numerical examples include the
time-dependent Stokes equations and vortex-dominated flows for small viscosities, and the Beltrami
flow. The tests serve to verify pressure robustness and the ability of the methods to
preserve structures. Furthermore, the 3D Taylor-Green vortex is utilized as a working horse to
demonstrate that with a pressure robust DG method, a robust method for underresolved turbulent
incompressible flows can be realized. At the end we conclude in Section \ref{sec:conclusion-outlook}.

\section{Discontinuous Galerkin spatial discretization}
\label{sec:disc-galerk-spat}

The instationary incompressible Navier-Stokes equations in an open and bounded domain
\(\Omega\subset\mathbb{R}^d\) (\(d=2,3\)) and time interval \((0,T]\) with velocity \(v\) and pressure \(p\) as
unknowns for given right-hand side \(f\), viscosity \(\mu\) and density \(\rho\) are given by
\begin{subequations}
\label{eq:incomp_navierstokes}
\begin{align}
  \rho\partial_t v - \mu\Delta v + \rho (v\cdot\nabla) v + \tnabla p &= f
  && \textup{in} \; \Omega \times (0,T] \\
  \nabla\cdot v &= 0 && \textup{in} \; \Omega \times (0,T]\\
  v &= v_0 && \textup{for} \; t=0 \; . \\
\intertext{Either Dirichlet boundary condition for the velocity:}
  v &= g  && \textup{on} \; \Gamma_D, t\in (0,T] \label{eq:Dirichlet1}\\
\intertext{and free-slip boundary condition in addition if \(\Gamma_D\neq\partial\Omega\):}
  v\cdot\normal = 0 \quad \textup{and} \quad \mu\nabla v\normal\times\normal &= 0
  && \textup{on} \; \Gamma_S = \partial\Omega\setminus\Gamma_D , t\in(0,T] \label{eq:freeslip1} \\
\intertext{together with}
  \int_\Omega p \dl x &= 0 && \textup{for all} \; t\in (0,T] \label{eq:Dirichlet2}\\
\intertext{or mixed boundary conditions:}
  v &= g && \textup{on} \; \Gamma_D\neq\{\partial\Omega, \emptyset\} \label{eq:Mixed1}\\
  v\cdot\normal = 0 \quad \textup{and} \quad \mu\nabla v\normal\times\normal &= 0
  && \textup{on} \; \Gamma_S \; \textup{disjoint to} \; \Gamma_D \label{eq:freeslip2} \\
  \mu\nabla v \normal - p\normal -\rho\beta(v\cdot\normal)_{-}v &= 0
  &&\textup{on} \; \Gamma_N = \partial\Omega\setminus (\Gamma_D\cup\Gamma_S) \label{eq:Mixed2}
\end{align}
\end{subequations}
are supplemented with the system.
For pure Dirichlet boundary conditions \(g\) is required to satisfy the compatibility condition
\(\int_{\partial\Omega} g\cdot\normal \dl s = 0\) and in addition of free-slip boundary conditions
\(\int_{\Gamma_D} g\cdot\normal \dl s = 0\). On \(\Gamma_N\), \((v\cdot \normal)_{-} = \min(0,v\cdot\normal)\) denotes the
negative part of the flux across the boundary. The parameter \(\beta\) can take the values \(\beta=0\)
(classical do-nothing, CDN) or \(\beta=\frac12\) (directional do-nothing, DDN), \cite{Braack2014}.
In the numerical examples below we will also consider periodic boundary
conditions in addition.  Under appropriate assumptions the Navier-Stokes
problem in weak form has a solution in $(H^1(\Omega))^d\times
L^2(\Omega)$ for $t\in(0,T]$, \cite{GiraultRaviartBook,TemamBook}. In absence of do-nothing conditions the
pressure is only determined up to a constant and is in the space $L^2_0(\Omega) = \{q\in L^2(\Omega) \mid \int_\Omega q \dl x = 0 \}$.

For the discretization let $\mathcal{E}_h$ be a quadrilateral mesh (in dimension \(d=2\)) or a
hexahedral mesh (in dimension \(d=3\)) with maximum diameter \(h\). We denote by $\Gamma_h^\interior$ the
set of all interior faces, by $\Gamma_h^{D}$ the set of all faces intersecting with the Dirichlet boundary
$\Gamma_D$, by \(\Gamma_h^S\) the set of all faces intersecting with the free-slip boundary and by $\Gamma_h^N$ the
set of all faces intersecting with the mixed boundary \(\Gamma_h^N\). We set $\Gamma_h = \Gamma_h^\interior \cup \Gamma_h^D \cup
\Gamma_h^S \cup \Gamma_h^N$. To an interior face $e\in\Gamma_h^\interior$ shared by elements $E_e^1$ and $E_e^2$ we define an
orientation by its unit normal vector $\normal_e$ pointing from $E_e^1$ to $E_e^2$. The jump and average
of a scalar-valued function $\phi$ on a face is then defined by
\begin{align}
  [\phi] &= \phi\mid_{E_e^1} - \phi\mid_{E_e^2} \quad \; = \phi^\interior -
  \phi^\exterior, \label{eq:jump_average} \\
  \{\phi\} &= \frac12 \phi\mid_{E_e^1} + \frac12 \phi\mid_{E_e^2} =
  \frac12 \phi^\interior + \frac12 \phi^\exterior \notag \; .
\end{align}
Note that the definition of jump and average can be extended in a natural way to vector and
matrix-valued functions.  If $e\in\partial\Omega$ then $\normal_e$ corresponds to the outer normal
vector $\normal$.  Below we make heavy use of the identities and notation, respectively:
\begin{align}
  [u v] &= [u] \{v\} + \{u\} [v] \;,
  & (u,v)_{0,\omega} &= \int_{\omega} u v\,\dl x\; ,
  &&\text{($u,v$ scalar-valued)} \label{eq:jump_of_product}\\
  [u \cdot v] &= [u] \cdot \{v\} + \{u\} \cdot [v] \;,
  & (u,v)_{0,\omega} &= \int_{\omega} u \cdot v\,\dl x\; ,
  &&\text{($u,v$ vector-valued)} \notag \\
  [u : v] &= [u] : \{v\} + \{u\} : [v] \;,
  & (u,v)_{0,\omega} &= \int_{\omega} u : v\,\dl x\; ,
  &&\text{($u,v$ matrix-valued)} \notag .
\end{align}
where \(\omega\subset\Omega\) is a \(d\)-dimensional subset together with the \(d\)-dimensional
measure \(\dl x\). The same shorthand notation holds for the hypersurface measure \(\dl s\) when
integrating over codimension one subsets as (parts of) the boundary or possible collection of
faces. The DG discretization on hexahedral meshes is based on the non-conforming finite element
space of polynomial degree $p$
\begin{align}
Q_h^p = \{ v\in L^2(\Omega) \mid v|_E= q\circ\mu_E^{-1},  q\in \mathbb{Q}_{p,d}, E\in\mathcal{E}_h \}
\end{align}
where $\mu_E : \hat E \to E$ is the transformation from the reference cube $\hat E$ to $E$
and $\mathbb{Q}_{p,d}$ is the set of polynomials of maximum degree $p$ in $d$ variables.
The approximation spaces for velocity and pressure are then
\begin{subequations}
\begin{align}
X^p_h \times M_h^{p-1} &= (Q_h^p)^d \times (Q_h^{p-1}\cap L^2_0(\Omega))
&&\text{(Dirichlet b. c.)},\\
X^p_h \times M_h^{p-1} &= (Q_h^p)^d \times Q_h^{p-1}
&&\text{(mixed b. c.)} .
\end{align}
\end{subequations}
We make use of the following mesh-dependent forms
defined on $X_h^p \times X_h^p$, $X_h^p \times M_h^{p-1}$, $X_h^p$ and $M_h^{p-1}$,
respectively:
\begin{subequations}
\begin{align}
  a(u,v) &= d(u,v) + J_0(u,v), \ \text{where}\\
  d(u,v) &= \sum_{E\in\mathcal{E}_h} (\mu\nabla u, \nabla v)_{0,E} \notag \\
  &\quad - \sum_{e\in\Gamma_h^\interior} (\mu\{\nabla u\}\normal_e, [v])_{0,e}
  - \sum_{e\in\Gamma_h^D} (\mu\nabla u^\interior\normal_e , v^\interior)_{0,e}
  - \sum_{e\in\Gamma_h^S} (\mu\normal_e^T\nabla u^\interior\normal_e , v^\interior\cdot\normal_e)_{0,e}, \label{eq:dg_discr_diffusion}\\
  J_0(u,v) &= \epsilon \sum_{e\in\Gamma_h^\interior} (\mu\{\nabla v\}\normal_e , [u])_{0,e}
  + \epsilon \sum_{e\in\Gamma_h^D} (\mu\nabla v^\interior\normal_e,u^\interior)_{0,e}
  + \epsilon \sum_{e\in\Gamma_h^S} (\mu\normal_e^T\nabla v^\interior\normal_e,u^\interior\cdot\normal_e)_{0,e} \notag\\
  &\quad + \sum_{e\in\Gamma_h^\interior} \mu\frac{\sigma}{h_e} ([u], [v])_{0,e}
  + \sum_{e\in\Gamma_h^D} \mu\frac{\sigma}{h_e} (u^\interior, v^\interior)_{0,e}
  + \sum_{e\in\Gamma_h^S} \mu\frac{\sigma}{h_e} (u^\interior\cdot\normal_e, v^\interior\cdot\normal_e)_{0,e}, \label{eq:dg_penalty} \\
  b(v,q) &= -\sum_{E\in\mathcal{E}_h} (\nabla\cdot v,q)_{0,E} +
  \sum_{e\in\Gamma_h^\interior} ([v]\cdot\normal_e,\{q\} )_{0,e}
  +\sum_{e\in\Gamma_h^D\cup\Gamma_h^S} (v^\interior \cdot \normal,q^\interior )_{0,e} ,\\
  l(v;t) &=\sum_{E\in\mathcal{E}_h} (f(t) , v)_{0,E}
  + \epsilon \sum_{e\in\Gamma_h^D} (\mu\nabla v^\interior\normal_e , g(t) )_{0,e}
  + \sum_{e\in\Gamma_h^D} \mu\frac{\sigma}{h_e} (g(t),v^\interior)_{0,e} ,\\
  r(q;t) &=\sum_{e\in\Gamma_h^D}  (g(t)\cdot\normal,q^\interior)_{0,e} \; .
\end{align}
\end{subequations}
Here we made the time dependence of the right hand side functionals explicit. For ease of writing
this will be omitted mostly below. In the interior penalty parameter \(\sigma / h_e\), the denominator
accounts for the mesh dependence. The formula for \(h_e\),
\begin{equation*}
  h_e = \begin{cases}
    \frac{\min\left( \left| E^{\interior}(e) \right|, \left|
          E^{\exterior}(e) \right| \right)}
    { \left| e \right| } &, E^{\interior}(e) \cap E^{\exterior}(e) = e \\
    \frac{ \left| E^{\interior}(e) \right| }{ \left| e \right| } &,
    E^{\interior}(e)
    \cap \Gamma_D = e
  \end{cases} \; ,
\end{equation*}
has been stated in \cite{HoustonHartmann2008} where it was proven that this choice ensures
coercivity of the bilinear form for anisotropic meshes.  For $\sigma$ we choose
\(\sigma = \alpha p(p+d-1)\) as in \cite{amg4dg} with \(\alpha\) a user-defined parameter to be
chosen \(\alpha=3\) in the computations reported below.  In $J_0$ the SIPG ($\epsilon = -1$) method is preferred
since the matrix of the linear system in absence of the convection term is then symmetric. Other
choices are the NIPG ($\epsilon = 1$) or IIPG ($\epsilon = 0$) method.

\cite{Hillewaert13} presents a rigorous analysis on the optimal penalty parameter where exact bounds
from the trace inverse inequality for triangles, tetrahedra, quadrilaterals, hexahedra, wedges and
pyramids are derived. See Table 3.1 in his doctoral dissertation. There are two conditions mentioned
to ensure coercivity, Equation (3.22) and Equation (3.23). \cite{Hillewaert13} also verifies that an
optimal penalty parameter is not sharply confined by these equations. The condition expressed in
(3.22) is used in \cite{KRANK2017634} and a related series of publications. The condition expressed
in (3.23) gives the same mesh dependence on the penalty parameter as the formula for \(h_e\). It is
cheaper than the former condition which requires to iterate over all faces in the adjacent elements
for each face. Further in Equation (3.23) setting the number of faces per element (\(n_e\) in his
notation where \(e\) denotes an element) equal to \(2d\) for quadrilaterals and hexahedra, and our
choice of \(\alpha = 3\), yields a good agreement compared to the penalty parameter \(\sigma/h_e\) introduced above.

The convective term in the Navier-Stokes equations can be written in conservative form as \(\nabla\cdot
F(v)\) with the convective flux matrix \(F(v) = v\otimes v\). In the DG method this term is then treated
with an upwind scheme introduced in \cite{2016arXiv161200657P}:
\begin{equation}
\begin{split}
c(v,\varphi) &= -\sum_{E\in\mathcal{E}_h} ( F(v) , \nabla\varphi)_{0,E} +
  \sum_{e\in\Gamma_h^{\text{int}}} (\hat{F}_e(v,\normal_e) , [\varphi ])_{0,e}
 + \sum_{e\in\Gamma_h^{D}\cup\Gamma_h^{N}}
 (\hat{F}_e(v,\normal_e),\varphi)_{0,e}\; ,
  \label{eq:part_integration_convection_conservative_final}
\end{split}
\end{equation}
with the numerical fluxes
\begin{equation*}
\hat F_e(v,\normal_e) = \left\{\begin{array}{ll}
\max(0,\{v\}\cdot\normal_e) v^{\text{int}} + \min(0,\{v\}\cdot\normal_e) v^\text{ext}  & e\in\Gamma_h^\text{int}\\
\max(0,v^\text{int}\cdot\normal_e) v^{\text{int}} + \min(0,v^\text{int}\cdot\normal_e) g & e\in\Gamma_h^D\\
\max(0,v^\text{int}\cdot\normal_e) v^{\text{int}} & e\in\Gamma_h^N
\end{array}\right. \; .
  \label{eq:vijayasundaram_flux}
\end{equation*}
On the outflow boundary the variational form of the DDN contribution is
\begin{equation}
  \label{eq:ddn_outflow}
  s_o(u,v) = \sum_{e\in\Gamma_h^N} ((u\cdot\normal)_{-}u,v)_{0,e} \; .
\end{equation}

The discrete in space, continuous in time formulation of the Navier-Stokes problem
(\ref{eq:incomp_navierstokes}) now seeks to find $v_h(t) : (0,T]\to X_h^p$,
$p_h(t) : (0,T]\to M_h^{p-1}$:
\begin{subequations}
\label{eq:discrete_problem}
\begin{align}
  \rho(\partial_t v_h, \varphi)_{0,\Omega} + a(v_h,\varphi) +
  \rho c(v_h,\varphi) - \rho\beta s_o(v_h,\varphi) + b(\varphi,p_h)
  &= l(\varphi;t) , \label{eq:variational_momentum} \\
  b(v_h,q) &= r(q;t) , \label{eq:variational_incomp}
\end{align}
\end{subequations}
for all $(\varphi,q)\in X_h^p\times M_h^{p-1}$. It can be shown that the scheme satisfies the local mass
conservation property
\begin{equation}
\sum_{e\in\Gamma_h^{\text{int}}\cap\partial E} (\{v \}\cdot n_e,1)_{0,e}
+ \sum_{e\in\Gamma_h^{N}\cap\partial E} (v\cdot \normal_e,1)_{0,e}
+ \sum_{e\in\Gamma_h^D\cap\partial E}  (g\cdot\normal,1)_{0,e} = 0 .
\label{eq:local_mass_conservation}
\end{equation}

\section{Splitting method}
\label{sec:splitting-method}

The splitting method for solving (\ref{eq:discrete_problem}) relies on the Helmholtz decomposition
of the velocity field to correct for the divergence constraint (\ref{eq:variational_incomp}). In this
section we introduce at first the notation needed for the description of the Navier-Stokes splitting
method. We summarize the entire fractional stepping, and towards the end of this section we review
some of the discrete Helmholtz decompositions available in the literature. The novel
\(H(\text{div})\) postprocessing technique is presented afterwards in Section
\ref{sec:helmh-flux-ravi}.

\subsection{Helmholtz decomposition}
\label{sec:helmh-decomp}

The Helmholtz decomposition states that any vector field in \(L^2(\Omega)^d\) can be decomposed into a
divergence-free contribution and an irrotational contribution. In order to define the decomposition,
boundary conditions on the pressure need to be enforced which are not part of the underlying
Navier-Stokes equations. Consider for simplicity \(\Gamma_N = \emptyset\), no mixed boundary conditions, and let
us denote the space of weakly divergence-free functions by
\begin{equation*}
  \label{eq:weakly_divfree_space}
  H(\Omega) := \{ v \in L^2(\Omega)^d \mid (v,\tnabla f)_{0,\Omega} - (g\cdot\normal,f)_{0,\Gamma_D} = 0 \; \forall f\in H^1(\Omega) \} \; .
\end{equation*}
In addition, we employ the pressure space
\begin{equation*}
  \label{eq:helmholtz_dirichlet_complement}
  \Psi_D(\Omega) := \{ q\in H^1(\Omega) \mid (q,1)_{0,\Omega} = 0 \} \; .
\end{equation*}
Then the fundamental theorem of vector calculus states that for any \(w\in L^2(\Omega)^d\) there are unique
functions \(v\in H(\Omega)\) and \(\psi\in \Psi_D(\Omega)\) such that \[ w = v + \tnabla\psi \; .\]
The irrotational contribution is computed by solving the variational problem
\begin{equation*}
  \label{eq:Poisson1}
  (\tnabla \psi,\tnabla q)_{0,\Omega} = (w,\tnabla q)_{0,\Omega} - (g\cdot\normal,q)_{0,\Gamma_D} \qquad \forall q\in \Psi_D(\Omega)
\end{equation*}
and it can be readily checked that \(v = w - \tnabla\psi\) is in the space \(H(\Omega)\). The equation for
\(\psi\) is the weak formulation of a Poisson equation with homogeneous Neumann boundary conditions on
\(\partial\Omega = \Gamma_h^D\cup\Gamma_h^S\). Note that in a pressure correction scheme the Helmholtz decomposition is
divided by the time step to obtain the physical pressure \(p\) in the Navier-Stokes system.

All the variants of discrete Helmholtz decomposition covered in this article are based on the
solution of a pressure Poisson equation, as in the continuous case. Therefore we define the form
\begin{equation}
  \label{eq:pressurepoisson_dg_laplacian}
  \begin{split}
    \alpha(\psi_h,q_h) &= \sum_{E\in\mathcal{E}_h} (\nabla\psi_h,\nabla q_h)_{0,E} - \sum_{e\in\Gamma_h^{\text{int}} } (\{\nabla\psi_h \}\cdot\normal_e,[q_h])_{0,e}
                 - \sum_{e\in\Gamma_h^{N}} (\nabla\psi_h\cdot \normal_e,q_h)_{0,e} \\
    &\hspace*{1.5em} - \sum_{e\in\Gamma_h^{\text{int}} } (\{\tnabla q_h \}\cdot \normal_e,[\psi_h])_{0,e} +
      \sum_{e\in\Gamma_h^{\text{int}} } \frac{\sigma}{h_e} ([\psi_h],[q_h])_{0,e} \\
    &\hspace*{1.5em} - \sum_{e\in\Gamma_h^{N}} (\tnabla q_h\cdot \normal_e,\psi_h)_{0,e} + \sum_{e\in\Gamma_h^{N}} \frac{\sigma}{h_e}
    (\psi_h,q_h)_{0,e} \; , \quad \psi_h,q_h\in M_h^{p-1} \; ,
  \end{split}
\end{equation}
with \(\sigma = \alpha(p-1)(p+d-2)\) the penalty parameter for the space \(M_h^{p-1}\).
This corresponds to the SIPG formulation, see \cite{ShuyuWheeler2005}, of Poisson's equation with homogeneous Neumann boundary
conditions on \(\Gamma_h^D\cup\Gamma_h^S\) and homogeneous Dirichlet boundary conditions on \(\Gamma_h^N\) (which might
be empty), c.f. \cite{2016arXiv161200657P}.

\subsection{Rotational Incremental Pressure Correction Scheme}
\label{sec:ripcs}

For ease of writing the DDN-term is omitted in the summary of the fractional step technique.
The incremental pressure correction scheme in rotational form
\cite{FLD:FLD373,2016arXiv161200657P} for solving (\ref{eq:discrete_problem}) reads:

\medskip
\noindent
Given \(v_h^k, p_h^k\) at time \(t^k\), compute \(v_h^{k+1}, p_h^{k+1}\)
at time \(t^{k+1} = t^k + \Delta t^{k+1}\) as follows:
\begin{enumerate}
\item Choose explicit extrapolation of pressure \(p_h^{\bigstar,k+1}\)
  at time \(t^{k+1}\) and compute tentative velocity $\tilde{v}_h^{k+1}$
  by temporal advancement of
  \begin{equation*}
    \rho(\partial_t v_h, \varphi_h)_{0,\Omega} + a(v_h,\varphi_h) +
    \rho c(v_h,\varphi_h) + b(\varphi_h,p_h^{\bigstar,k+1}) = l(\varphi_h;t) \quad \forall \varphi_h\in X_h^p \; .
  \end{equation*}
  with either Alexander's second order strongly S-stable scheme \cite{Alexander:1977:DIR} or the
  Fractional-Step \(\theta\)-method.
\item Perform discrete Helmholtz decomposition to obtain \(v_h^{k+1} = \mathcal{P}_h\tilde v_h^{k+1}\)
  and pressure correction \(\delta p_h^{k+1}\) using Algorithm \ref{alg:divdiv_proj},
  \ref{alg:divdivconti_proj}, \ref{alg:RT_reconstruction}
  or Algorithm \ref{alg:HelmholtzRT_reconstruction}, to be described below.

  This involves the solution of a pressure Poisson equation with unknown \(\delta p_h^{k+1} = \psi_h / \Delta t^{k+1}\)\,:
  \begin{align*}
    \alpha(\delta p_h^{k+1}, q_h)
    &= \frac{1}{\Delta t^{k+1}} (b(\tilde{v}_h^{k+1}, q_h) - r(q_h;t^{k+1}))
    \quad \forall q_h\in M_h^{p-1},
  \end{align*}
  and a postprocessing step.
\item Update to new pressure \(p_h^{k+1}\):
  \begin{equation*}
    (p_h^{k+1}, q_h) = (\omega \delta p_h^{k+1} + p_h^{\bigstar,k+1}, q_h)
    + \mu (b(\tilde{v}_h^{k+1}, q_h) - r(q_h;t^{k+1})) \quad \forall q_h\in M_h^{p-1} \; .
  \end{equation*}
\end{enumerate}
For the second order formulation of the splitting method, we set \(\omega = \frac32\) and
\(p_h^{\bigstar,k+1} = p_h^k\). The last term in the pressure update is specific for the rotational
form of the incremental pressure correction scheme. It ensures that the splitting reproduces
stationary solutions of the Navier-Stokes equations.

\subsection{Stabilization enhanced projections}
\label{sec:stab-enhanc-proj}

The authors of \cite{KRANK2017634} have proposed discrete Helmholtz decompositions using
stabilization terms. One such variant adds the same term to the solenoidal projection step as
grad-div stabilization to the momentum equation whereupon the latter is commonly used within coupled
solution approaches in order to weakly enforce exactly divergence-free solutions. We refer to this
variant as div-div projection. Another variant further adds a normal continuity penalty term to the
solenoidal projection step that weakly enforces \(H(\text{div})\)-conformity. We refer to this
variant as div-div-conti projection.
\begin{center}
\captionsetup{style=ruled,type=algorithm,skip=0pt}
\makeatletter
  \fst@algorithm\@fs@pre
\makeatother
\caption{
  The div-div projection is given by the following algorithm:
}
\makeatletter
  \@fs@mid\vspace{2pt}
\makeatother
\label{alg:divdiv_proj}
\begin{enumerate}[i)]
\item For any tentative velocity \(w_h\in X_h^p\) and fixed \(t\) solve
\begin{equation*}
\psi_h\in M_{h}^{p-1}: \quad \alpha(\psi_h,q_h)= b(w_h,q_h)-r(q;t) \quad \forall q_h\in M_{h}^{p-1}.
\end{equation*}
\item Set \(\mathcal{P}_h w_h = v_h\) where \(v_h\) solves
\begin{equation}
\label{eq:divdiv_subtraction}
(v_h,\varphi_h)_{0,\Omega} + \sum_{E\in\mathcal{E}_h}\tau_{D,E} (\nabla\cdot v_h,\nabla\cdot\varphi_h)_{0,E}
= (w_h,\varphi_h)_{0,\Omega} - (\tnabla \psi_h,\varphi_h)_{0,\Omega}
\quad \forall \varphi_h\in X_h^p
\end{equation}
where \(\tau_{D,E}\) is a per-cell penalization constant.

This requires the solution of an element-local system which is not diagonal. As reported in
\cite{KRANK2017634,FEHN2018667,2016arXiv161200657P} this gives good results with quite small
pointwise divergence. However, the projected velocity does not satisfy a local mass conservation
property and \((\mathcal{P}_h)^2\neq \mathcal{P}_h\).
\end{enumerate}
\makeatletter
  \vspace{2pt}\@fs@post
\makeatother
\end{center}
\begin{center}
\captionsetup{style=ruled,type=algorithm,skip=0pt}
\makeatletter
  \fst@algorithm\@fs@pre
\makeatother
\caption{
  The div-div-conti projection is given by the following algorithm:
}
\makeatletter
  \@fs@mid\vspace{2pt}
\makeatother
\label{alg:divdivconti_proj}
\begin{enumerate}[i)]
\item For any tentative velocity \(w_h\in X_h^p\) and fixed \(t\) solve (same as before)
\begin{equation*}
\psi_h\in M_{h}^{p-1}: \quad \alpha(\psi_h,q_h)= b(w_h,q_h)-r(q;t) \quad \forall q_h\in M_{h}^{p-1}.
\end{equation*}
\item Set \(\mathcal{P}_h w_h = v_h\) where \(v_h\) solves
\begin{equation}
\label{eq:divdivconti_subtraction}
\begin{split}
(v_h,\varphi_h)_{0,\Omega} +& \sum_{E\in\mathcal{E}_h}\tau_{D,E} (\nabla\cdot v_h,\nabla\cdot\varphi_h)_{0,E} + \sum_{e\in\Gamma_h^\interior}\tau_{C,e}([v_h]\cdot\normal_e,[\varphi_h]\cdot\normal_e)_{0,e} \\
&\hspace*{1.5em} = (w_h,\varphi_h)_{0,\Omega} - (\tnabla \psi_h,\varphi_h)_{0,\Omega} \quad \forall \varphi_h\in X_h^p
\end{split}
\end{equation}
where \(\tau_{D,E}\) is a per-cell penalization and \(\tau_{C,e}\) is a per interior face
penalization constant.

This gives good results with small pointwise divergence as well. The authors of \cite{FEHN2018667}
further demonstrate that a robust DG method for underresolved turbulent flow can be
realized. However, the projected velocity does not satisfy a local mass conservation property and
\((\mathcal{P}_h)^2\neq \mathcal{P}_h\).
\end{enumerate}
\makeatletter
  \vspace{2pt}\@fs@post
\makeatother
\end{center}
\noindent
The continuity penalty introduces inter-element couplings such that the system is not
block-diagonal. This system has the same stencil as a Poisson operator and the convergence of the
CG solver employed with an iterative block Jacobi preconditioner \cite{BASTIAN2019417} does
deteriorate. We will mostly utilize the div-div projection in this work.

\subsection{Pressure Poisson Raviart-Thomas projection}
\label{sec:press-poiss-ravi}

Another variant has been proposed by the authors \cite{2016arXiv161200657P}. It reconstructs the
negative gradient of the pressure correction in the Raviart-Thomas space of degree \(k\) before
subtracting the irrotational contribution. On hexahedral meshes these spaces are given by
\cite{Brezzi91}:
\begin{equation}
\RT_h^k = \{ v\in H(\text{div};\Omega) \mid v|_E\in \RT_E^k \, \forall E\in\mathcal{E}_h\}
\end{equation}
with the Raviart-Thomas space on element \(E\) given by
\begin{equation}\label{eq:RTonE}
\RT_E^k = \{ v\in H(\text{div};E) \mid v = T_E(\hat v),
(\hat v)_i = \sum_{\{ \alpha \mid 0\leq\alpha_j\leq k+\delta_{ij}\}} c_{i,\alpha} \hat x^\alpha\}
\end{equation}
where we made use of the Piola transformation to the element \(E\in\mathcal{E}_h\),
i.e. for \(\mu_E(\hat x): \hat E \rightarrow E\) defined as
\begin{equation*}
T_E(\hat v)(x) = \frac{1}{|\det \nabla\mu_E(\hat x)|} \nabla\mu_E(\hat x) \hat v(\hat x) \; .
\end{equation*}
For \(k>0\) the construction needs also the space
\begin{equation}
\Psi_E^k = \{ v\in H(\text{div};E) \mid v = T_E(\hat v),
(\hat v)_i = \sum_{\{ \alpha \mid 0\leq\alpha_j\leq k-\delta_{ij}\}} c_{i,\alpha} \hat x^\alpha\} \; .
\end{equation}
Note that in contrast to \eqref{eq:RTonE} the polynomial degree in direction \(i\) in component
\(i\) is decreased instead of increased.

Assume that \(\psi_h\in M_h^{p-1}\) solves the DG-discretized pressure Poisson equation. Following
\cite{ErnHdiv2007} we now compute \(\gamma_h = G_h\psi_h\in\RT_h^k\) as reconstruction of \(-\tnabla \psi_h\) as
follows. On element \(E\in\mathcal{E}_h\) with faces \(e\in\partial E\) define
\begin{subequations}
\label{eq:RT_reconstruction}
\begin{align}
(\gamma_h\cdot \normal_e,q)_{0,e} &=
(-\{ \nabla\psi_h \}\cdot \normal_e + \frac{\sigma}{h_e} [\psi_h],q)_{0,e} && e\in\Gamma_h^\text{int}, q\in Q^{k}_e,\\
(\gamma_h\cdot \normal_e,q)_{0,e} &=
(- \nabla\psi_h \cdot \normal_e + \frac{\sigma}{h_e} \psi_h,q)_{0,e} && e\in\Gamma_h^N, q\in Q^{k}_e,\\
(\gamma_h\cdot \normal_e,q)_{0,e} &= 0 && e\in\Gamma_h^D\cup\Gamma_h^S, q\in Q^{k}_e, \label{eq:rt_bc}
\end{align}
and for $k>0$ define in addition
\begin{equation}
\begin{split}
(\gamma_h,r)_{0,E} &= -( \nabla\psi_h,r)_{0,E}
+\frac12 \sum_{e\in\partial E \cap \Gamma_h^\text{int}} (r\cdot \normal_e,[ \psi_h ])_{0,e}
 +\sum_{e\in\partial E \cap \Gamma_h^N} (r\cdot \normal_e, \psi_h)_{0,e} ,
\hspace*{0.75em} \forall r \in \Psi_E^k \; .
\label{eq:RT_projection_c}
\end{split}
\end{equation}
\end{subequations}
We can now define this discrete Helmholtz decomposition:
\begin{center}
\captionsetup{style=ruled,type=algorithm,skip=0pt}
\makeatletter
  \fst@algorithm\@fs@pre
\makeatother
\caption{
  The pressure Poisson Raviart-Thomas projection is given by the following algorithm:
}
\makeatletter
  \@fs@mid\vspace{2pt}
\makeatother
\label{alg:RT_reconstruction}
\begin{enumerate}[i)]
\item For any tentative velocity \(w_h\in X_h^p\) and fixed \(t\) solve (same as before)
\begin{equation*}
\psi_h\in M_{h}^{p-1}: \quad \alpha(\psi_h,q_h)= b(w_h,q_h)-r(q_h;t) \quad \forall q_h\in M_{h}^{p-1}.
\end{equation*}
\item Reconstruct \(\gamma_h = G_h\psi_h \in \RT_h^k\).
\item Set \(\mathcal{P}_h w_h = v_h\) where \(v_h\) solves
\begin{equation*}
(v_h,\varphi_h)_{0,\Omega}  = (w_h,\varphi_h)_{0,\Omega} + (G_h\psi_h,\varphi_h)_{0,\Omega}
\quad \forall \varphi_h\in X_h^p .
\end{equation*}
This requires the solution of a (block-) diagonal system.
\end{enumerate}
\makeatletter
  \vspace{2pt}\@fs@post
\makeatother
\end{center}
\noindent
In \cite{2016arXiv161200657P} it was shown on rectangular quadrilateral/hexahedral meshes that if the
flux is reconstructed in \(\RT_h^k,\; k=p-1\), \(v_h = \mathcal{P}_h w_h\) satisfies the discrete
continuity equation exactly, i.e.
\begin{equation*}
  b(\mathcal{P}_h w_h,q_h) = r(q_h;t) \quad\forall q_h\in M_h^{p-1} \; ,
\end{equation*}
and the operator \(\mathcal{P}_h\) therefore defines a projection. However, the pointwise divergence
of \(v_h\) is not zero, with values larger than obtained by div-div projection
(c.f. \cite{2016arXiv161200657P}). Note that the variant to be introduced in Section
\ref{sec:helmh-flux-ravi} successfully remedies on this point. In terms of local mass conservation,
it has been observed that it is sufficient to reconstruct the negative pressure gradient in the
Raviart-Thomas space of degree \(k=p-2\), and thus we will consider \(\RT_h^{p-2}\) for the pressure
Poisson flux in this work.

\section{Helmholtz flux Raviart-Thomas projection}
\label{sec:helmh-flux-ravi}

In this section we present a reconstruction of the Helmholtz flux \(w_h - \tnabla_h\psi_h\) in the
Raviart-Thomas of degree \(p-1\) that does not only satisfy the discrete continuity equation but is
also pointwise divergence-free. The development of Helmholtz flux Raviart-Thomas projection
originated from studying pressure robust discretizations of the incompressible Navier-Stokes
equations. In fact, as the numerical experiments show, the splitting scheme obtained by
Helmholtz flux Raviart-Thomas projection evidences to be pressure robust. We start the derivation by
considering a reconstruction operator that is divergence-preserving. With the help of a
divergence-preserving reconstruction operator, recently, the following discretizations were shown to
give pressure robustness: \cite{LINKE2014782} presents a modified Crouzeix-Raviart element for the
incompressible Navier-Stokes equations. The authors of \cite{doi:10.1137/15M1047696} derive a
reconstruction operator for the mixed finite element pair on triangles consisting of conforming
\(\mathcal{P}_2\) space for velocity enriched with bubble functions and discontinuous
\(\mathcal{P}_1\) for pressure. In \cite{DIPIETRO2016175} a higher order reconstruction operator for
a discontinuous method is presented on simplicial meshes where the cell-based unknowns are
eliminated by static condensation. The authors of \cite{doi:10.1137/16M1089964} develop a
reconstruction for pressure robust Stokes discretizations with continuous pressure finite
elements. However, to the best of our knowledge, no results have been presented on
quadrilateral/hexahedral meshes for discontinuous velocity and pressure spaces.

In order to present a divergence-preserving operator for both discontinuous velocity and pressure
spaces, let us recapitulate the notion of the discrete divergence operator \(B_h\). It is a map
\(B_h: X_h^p\rightarrow Q_h^{p-1}\) satisfying
\begin{equation}
  (B_h w_h, q_h)_{0,\Omega} = -b(w_h,q_h) + r(q_h) \quad \forall q_h\in Q_h^{p-1} \; .
  \label{eq:discrete_divergence_operator}
\end{equation}
The kernel of \(B_h\) is denoted by \(X_{h,div}^p\) and called the set of discretely divergence-free
vector fields.
Next, we introduce the reconstruction operator \(\Pi_h^{\RT_h^{p-1}}\) that maps elements of the
velocity space \(X_h^p\) to the Raviart-Thomas space \(\RT_h^{p-1}\) of degree \(p-1\). Based on the
results presented in \cite{LINKE2014782,doi:10.1137/15M1047696,DIPIETRO2016175}, we define
\(v_h = \Pi_h^{\RT_h^{p-1}}(w_h), \, w_h\in X_h^p,\) as follows: On element \(E\in\mathcal{E}_h\) with
faces \(e\in\partial E\) compute
\begin{subequations}
\begin{align}
  (v_h\cdot n_e,q)_{0,e} &= (\{w_h\}\cdot n_e,q)_{0,e} && e\in\Gamma_h^\text{int}, q\in Q^{p-1}_e \; , \label{eq:velocityrt_interface} \\
  (v_h\cdot n_e,q)_{0,e} &= (g\cdot n_e,q)_{0,e} && e\in\Gamma_h^D, q\in Q^{p-1}_e \; , \\
  (v_h\cdot n_e,q)_{0,e} &= 0 && e\in\Gamma_h^S, q\in Q^{p-1}_e \; , \label{eq:velocityrt_bcfreeslip} \\
  (v_h\cdot n_e,q)_{0,e} &= (w_h\cdot n_e,q)_{0,e} && e\in\Gamma_h^N, q\in Q^{p-1}_e \\
  \intertext{and in addition for \(p-1>0\)}
  (v_h,r)_{0,E} &= (w_h,r)_{0,E} && \forall r \in \Psi_E^{p-1} \; . \label{eq:velocityrt_volume}
\end{align}
\end{subequations}
Obviously this defines a projection. On rectangular quadrilateral/hexahedral meshes we are able to show:
\begin{theorem}
  \label{theo:divergence_preserving_operator}
  The velocity reconstruction operator \(\Pi_h^{\RT_h^{p-1}}\) defined by the equations
  (\ref{eq:velocityrt_interface})-~(\ref{eq:velocityrt_volume}) is divergence-preserving.
  In explicit it holds
  \begin{equation}
    \label{eq:divergence_preserving}
    \nabla\cdot\Pi_h^{\RT_h^{p-1}} = B_h \; .
  \end{equation}
  \begin{proof}
    Let \(v_h = \Pi_h^{\RT_h^{p-1}}(w_h)\) for \(w_h\in X_h^p\). Essential steps in the proof are that
    for any \(q_h\in Q_h^{p-1} \Rightarrow \tnabla q_h|_E\in\Psi_E^{p-1}\), and the alternative definition of the
    form \(b\) using integration by parts.
    \begin{align*}
      (\nabla\cdot v_h,q_h)_{0,\Omega} &= \sum_{E\in\mathcal{E}_h} (\nabla\cdot v_h,q_h)_{0,E} \\
      &= \sum_{E\in\mathcal{E}_h} -(v_h,\tnabla q_h)_{0,E} + \sum_{e\in\partial E} (v_h\cdot n_e,q_h)_{0,e} \\
      &= -\sum_{E\in\mathcal{E}_h} (v_h,\tnabla q_h)_{0,E} + \sum_{e\in\Gamma_h} (v_h\cdot n_e,[q_h])_{0,e} \\
      &\hspace*{-0.4em}\stackrel{~(\ref{eq:velocityrt_volume})}= -\sum_{E\in\mathcal{E}_h} (w_h,\tnabla q_h)_{0,E} + \sum_{e\in\Gamma_h^N} (w_h\cdot n_e,q_h)_{0,e} \\
      &\qquad + \sum_{e\in\Gamma_h^\interior} (\{w_h\}\cdot n_e,[q_h])_{0,e} + \sum_{e\in\Gamma_h^D} (g(t)\cdot n_e,q_h)_{0,e} \\
      &= -b(w_h,q_h) + r(q_h) \\
      &= (B_hw_h,q_h)_{0,\Omega} \; .
    \end{align*}
    Since \(\nabla\cdot v_h\in Q_h^{p-1}\), we finally obtain \(\nabla\cdot\Pi_h^{\RT_h^{p-1}} = B_h\).
  \end{proof}
\end{theorem}
\begin{lemma}
  \label{lemma:discrete_divfree_to_divfree}
  The operator \(\Pi_h^{\RT_h^{p-1}}\) maps discretely divergence-free vector fields to
  divergence-free vector fields. The image also satisfies the discrete continuity equation.
  \begin{proof}
    Let \(v_h = \Pi_h^{\RT_h^{p-1}}(w_h)\) for \(w_h\in X_{h,div}^p\). Then according to Theorem
    \ref{theo:divergence_preserving_operator} \(\nabla\cdot v_h = B_h w_h = 0\). To show that \(v_h\)
    satisfies the discrete continuity equation, we use \([v_h]\cdot\normal_e=0\) on \(e\in\Gamma_h^\interior\)
    and (\ref{eq:velocityrt_bcfreeslip}) to add zeroes
    \begin{align*}
      0 &= (\nabla\cdot v_h,q_h)_{0,\Omega} \\
        &= (\nabla\cdot v_h,q_h) - \sum_{e\in\Gamma_h^\interior} ([v_h]\cdot\normal_e,\{q_h\})_{0,e}
          - \sum_{e\in\Gamma_h^S} (v_h\cdot\normal_e,q_h)_{0,e} \\
        &\hspace*{1.5em} - \sum_{e\in\Gamma_h^D} (v_h\cdot\normal_e,q_h)_{0,e} + \sum_{e\in\Gamma_h^D} (v_h\cdot\normal_e,q_h)_{0,e} \\
        &= -b(v_h,q_h) + r(q_h) \; .
    \end{align*}
  \end{proof}
\end{lemma}

The key idea for Helmholtz flux Raviart-Thomas projection originates from improving the
divergence of \(v_h = w_h + G_h\psi_h\), the velocity field obtained by pressure Poisson Raviart-Thomas
projection. In \cite{2016arXiv161200657P}, Lemma 2, it was shown that the divergence does not
vanish pointwise but is controlled in an integral sense by the jumps of the tentative
velocity. Based on the discussion above, a straightforward choice is to take the image of the
divergence-preserving reconstruction operator. This will lead to a velocity field that is pointwise
divergence-free. Hence, the proposition for Helmholtz flux Raviart-Thomas projection consists of
combining the reconstruction operator for the velocity field with the accurate pressure Poisson flux
reconstruction: \(v_h = \Pi_h^{\RT_h^{p-1}}(w_h) + G_h\psi_h\).

\noindent
Explicitly: On element \(E\in\mathcal{E}_h\) with faces \(e\in\partial E\) compute
\begin{subequations}
\begin{align}
  (v_h\cdot n_e,q)_{0,e} &= (\{w_h\}\cdot n_e,q)_{0,e} + (-\{\nabla\psi_h\}\cdot n_e + \frac{\sigma}{h_e}[\psi_h],q)_{0,e} && e\in\Gamma_h^\text{int}, q\in Q^{p-1}_e
  \; , \label{eq:helmholtzrt_interface} \\
  (v_h\cdot n_e,q)_{0,e} &= (g\cdot n_e,q)_{0,e} && e\in\Gamma_h^D, q\in Q^{p-1}_e \; , \\
  (v_h\cdot n_e,q)_{0,e} &= 0 && e\in\Gamma_h^S, q\in Q^{p-1}_e \; , \\
  (v_h\cdot n_e,q)_{0,e} &= (w_h\cdot n_e,q)_{0,e} + (-\nabla\psi_h\cdot n_e + \frac{\sigma}{h_e}\psi_h,q)_{0,e} && e\in\Gamma_h^N, q\in Q^{p-1}_e \\
  \intertext{and in addition for \(p-1>0\)}
  (v_h,r)_{0,E} &= (w_h,r)_{0,E} -( \nabla\psi_h,r)_{0,E} \notag \\
  &\qquad +\frac12 \sum_{e\in\partial E \cap \Gamma_h^\text{int}} (r\cdot \normal_e,[ \psi_h ])_{0,e}
  +\sum_{e\in\partial E \cap \Gamma_h^N} (r\cdot \normal_e, \psi_h)_{0,e}  && \forall r \in \Psi_E^{p-1} \; . \label{eq:helmholtzrt_volume}
\end{align}
\end{subequations}
Then the discrete Helmholtz decomposition is defined as:
\begin{center}
  \captionsetup{style=ruled,type=algorithm,skip=0pt}
  \makeatletter
    \fst@algorithm\@fs@pre
  \makeatother
  \caption{
    The Helmholtz flux Raviart-Thomas projection is given by the following algorithm:
  }
  \makeatletter
    \@fs@mid\vspace{2pt}
  \makeatother
  \label{alg:HelmholtzRT_reconstruction}
  \begin{enumerate}[i)]
  \item For any tentative velocity \(w_h\in X_h^p\) and fixed \(t\) solve (same as above)
    \begin{equation*}
      \psi_h\in M_h^{p-1}: \quad \alpha(\psi_h,q_h) = b(w_h,q_h) - r(q_h;t) \quad \forall q_h\in M_h^{p-1} \; .
    \end{equation*}
  \item Reconstruct \(\Pi_h^{\RT_h^{p-1}}(w_h) + G_h\psi_h\in\RT_h^{p-1}\).
  \item Set \(\mathcal{P}_hw_h = v_h\) where \(v_h\) solves
    \begin{equation*}
      (v_h,\varphi_h)_{0,\Omega} = \big(\Pi_h^{\RT_h^{p-1}}(w_h) + G_h\psi_h, \varphi_h\big)_{0,\Omega} \quad \forall\varphi_h\in X_h^p \; .
    \end{equation*}
    This requires the solution of a (block-) diagonal system.
  \end{enumerate}
  \makeatletter
    \vspace{2pt}\@fs@post
  \makeatother
\end{center}
\noindent
We are now able to prove the following theorem given on rectangular quadrilateral/hexahedral meshes.
\begin{theorem}
  \label{theo:helmholtzrt_projection}
  The operator \(\mathcal{P}_h\) introduced in Algorithm \ref{alg:HelmholtzRT_reconstruction} is a
  projection. The image \(v_h = \mathcal{P}_h w_h\) is pointwise divergence-free and satisfies the
  discrete continuity equation exactly.
  \begin{proof}
    First consider the divergence of \(v_h\). Owing to \(\nabla\cdot (G_h\psi_h) = -B_h w_h\) (shown in \cite{2016arXiv161200657P})
    and to Theorem \ref{theo:divergence_preserving_operator}
    \begin{align*}
      (\nabla\cdot v_h,q_h)_{0,\Omega} &= \big(\nabla\cdot\Pi_h^{\RT_h^{p-1}}(w_h) + \nabla\cdot (G_h\psi_h),q_h\big)_{0,\Omega} \\
      &= -b(w_h,q_h) + r(q_h) + b(w_h,q_h) - r(q_h) \\
      &= 0 \quad \forall q_h\in Q_h^{p-1} \; .
    \end{align*}
    Using this and following the line of argument in Lemma \ref{lemma:discrete_divfree_to_divfree}, we can
    further conclude
    \begin{alignat*}{3}
      b(v_h,q_h) &= b \big(\Pi_h^{\RT_h^{p-1}}(w_h),q_h\big) &&+ b(G_h\psi_h,q_h) \\
      &= \sum_{e\in\Gamma_h^D} (g\cdot n_e,q_h)_{0,e} &&+ 0 &&= r(q_h) \quad \forall q_h\in Q_h^{p-1} \; .
    \end{alignat*}
    A second application of \(\mathcal{P}_h\) produces at first a zero Helmholtz
    correction and leaves the input velocity unchanged since \(\Pi_h^{\RT_h^{p-1}}\) is a
    projection.
  \end{proof}
\end{theorem}

In contrast to the sole pressure Poisson flux reconstruction, those conservation properties cannot
be achieved when reconstructing in a Raviart-Thomas space of degree \(k<p-1\). For lower degrees
accuracy and especially divergence-preservation are lost. In the following two sections we will
illustrate conservation properties and temporal convergence of this splitting scheme.

\subsection{Numerical conservation properties}
\label{sec:numer-cons-prop}

We redo the numerical experiments on local mass conservation that have been carried out in
\cite{2016arXiv161200657P}, here for the novel Helmholtz flux Raviart-Thomas reconstruction. We
consider again the instationary Navier-Stokes equations on the domain \(\Omega = (-1,1)^2\) and take the
vortex decay given by the analytical solution of the 2D Taylor-Green vortex. We set \(\nu = 1/100\)
and do computations, as before, on a \(160\times 160\) rectangular mesh.

\begin{figure}[!htb]
  \includegraphics[width=0.5\textwidth]{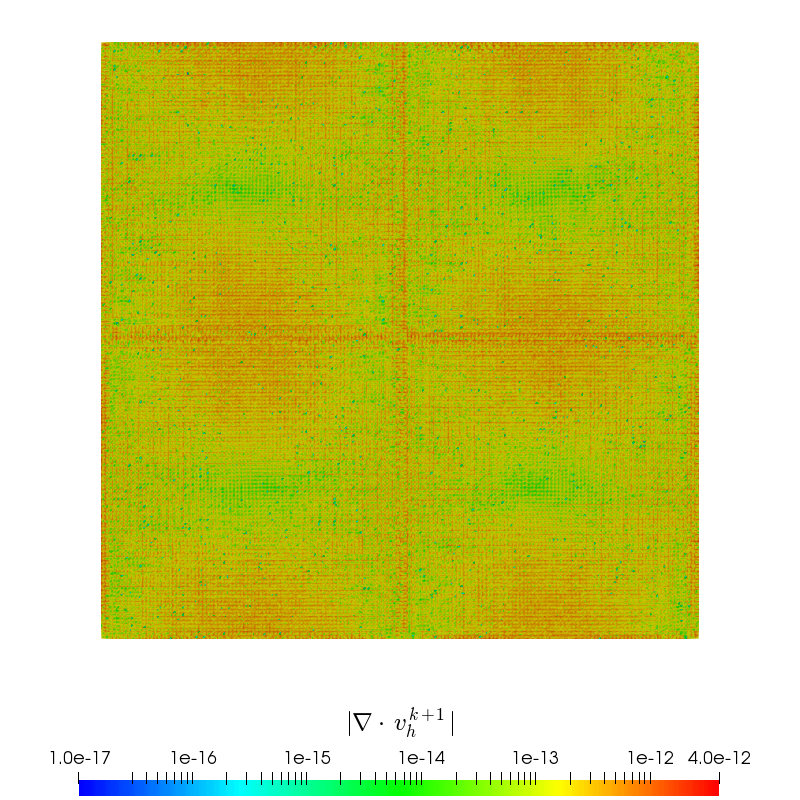}
  \includegraphics[width=0.5\textwidth]{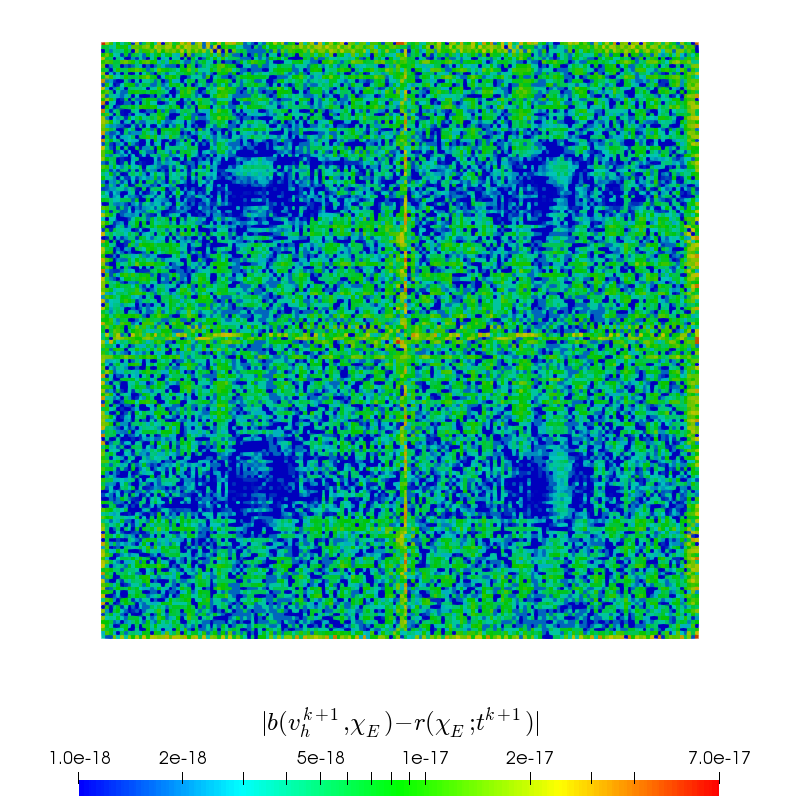}
  \caption{Pointwise divergence (left) and local mass conservation (right) of the 2D Taylor-Green
    vortex solution at time \(t = 1\), \(\Delta t = 0.025\) with the \THcube{2}{1} discretization and
    Helmholtz flux reconstruction in \(\RT_h^1\).
    \newline
    For comparison, the results of the other two discrete Helmholtz decompositions can be found in
    \cite{2016arXiv161200657P}.
  }
  \label{fig:conservation_properties_helmholtzflux_hdiv}
\end{figure}

As a comparison to the other discrete Helmholtz decompositions, figure
\ref{fig:conservation_properties_helmholtzflux_hdiv} shows the pointwise divergence and the local
mass conservation error on each mesh element for \(p=2\), the same time snapshot and same time step
size. Note that the error on local mass conservation is considered here as the right-hand side of
(\ref{eq:local_mass_conservation}). We observe that the pointwise divergence is lower than in the
case of div-div projection, and therefore lower as well than in the case of pressure Poisson
Raviart-Thomas reconstruction. The error on local mass conservation, however, has the same order of
magnitude as pressure Poisson \(H(\text{div})\) postprocessing, and is therefore lower than div-div
projection.

\subsection{Temporal convergence within the RIPCS}
\label{sec:temp-conv-with-ripcs}

We repeat the numerical convergence tests on the RIPCS that have been carried out in
\cite{2016arXiv161200657P}, here for the novel Helmholtz flux Raviart-Thomas reconstruction as
well. As already indicated by the previous results, the temporal error dominates in the range of
those measurements, and hence the errors do not change significantly and convergence rates in time
are retained. Thus, we only present an updated analysis with the Beltrami flow in three dimensions
using Helmholtz flux reconstruction. For the description of Beltrami flow refer to Section
\ref{sec:beltrami-flow}. Again we set \(\rho = \mu = 1\) and run computations on \(50^3\) cubic mesh.

\begin{figure}[!htb]
  \centering
  \includegraphics[width=0.6\textwidth]{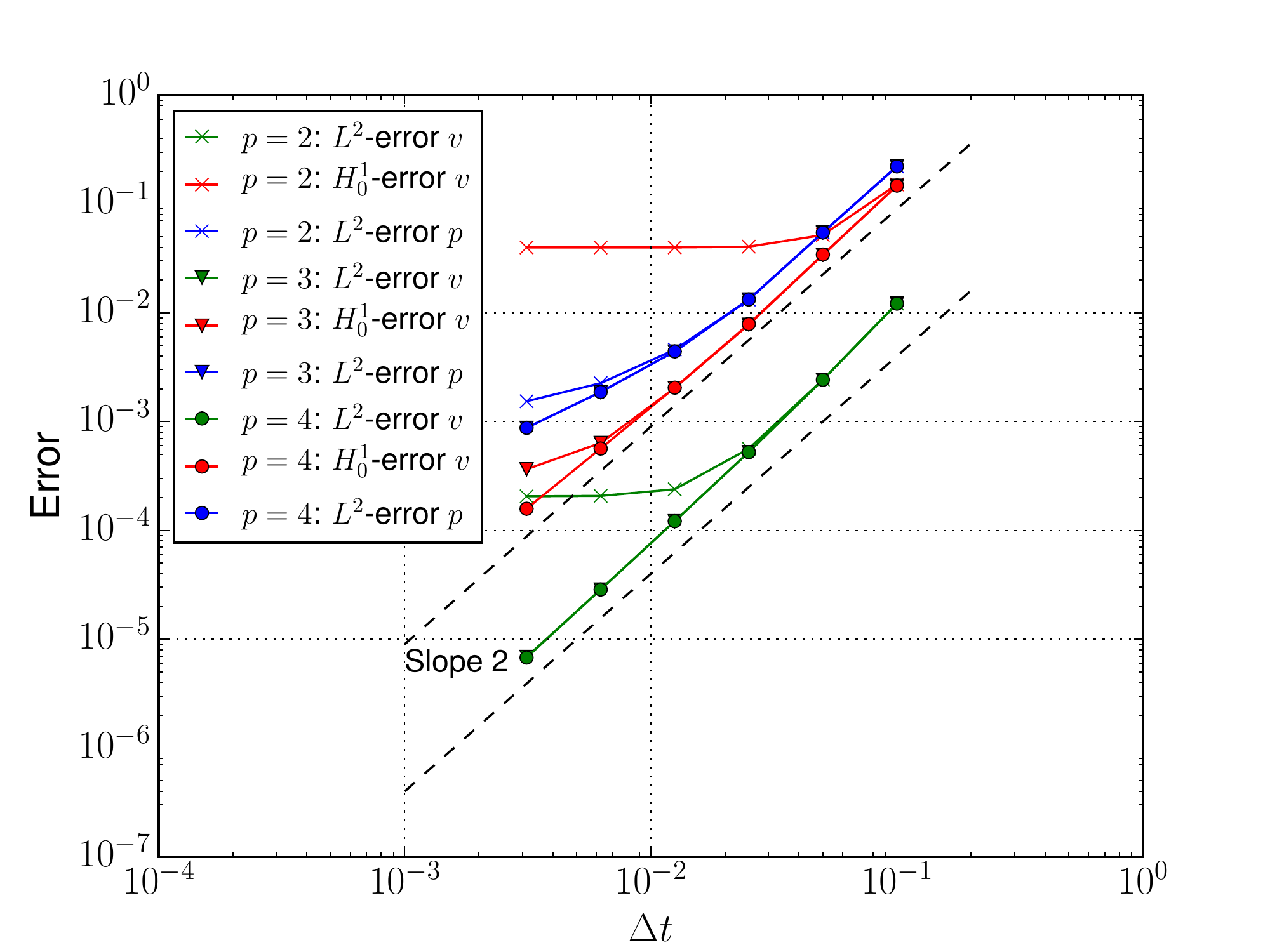}
  \caption{Errors and convergence rates at final time \(T=0.5\) for the Beltrami test problem using
    the spatial discretizations \THcube{p}{p-1}, \(p\in\{2,3,4\}\), and Helmholtz flux
    reconstruction in \(\RT_h^{p-1}\).
    \newline
    Previously obtained results with div-div projection can be found in \cite{2016arXiv161200657P}.}
  \label{fig:ripcs_beltrami_rates}
\end{figure}

Figure \ref{fig:ripcs_beltrami_rates} shows the errors and convergence rates as a function of
\(\Delta t\). The green curves show the \(L^2\)-error for the velocity, the red curves the \(H_0^1\)-error
for the velocity and the blue curves the \(L^2\)-error for the pressure obtained by the polynomial
degrees \(p=2,3,4\). Note that for \(p=2\) the spatial error in the splitting scheme becomes
all-dominant. This emerges at first for the \(H_0^1\)-error and \(L^2\)-error on the velocity. For
polynomial degree 4, however, the spatial error is negligible and the two error norms on the
velocity are perfectly second order convergent in time. The convergence for the pressure in
\(L^2\)-norm is close to \(\mathcal{O}(\Delta t^2)\).

\section{Numerical results}
\label{sec:numerical-results}

\subsection{Implementation}
\label{sec:implementation}

The parallel solver has been implemented in a high-performance C++ code using the object-oriented
DUNE finite element framework \cite{dune08:1,Bastian2016}. The incompressible Navier-Stokes solver
described uses the spectral discontinuous Galerkin method on quadrilateral/hexahedral meshes where
sum-factorization can be most easily applied to. The complexity reduction in the computation coming
from sum-factorization leads to a significant speedup especially for high-order DG
methods. Matrix-free methods are thereby put in a superior position to traditional matrix-based
methods as matrix entries can be computed faster on-the-fly than loaded from memory. For the
representation of function spaces we use tensor product bases built from Gauss-Lobatto-Lagrange
polynomials, for the evaluation of integrals we use (non-collocated) Gauss-Legendre quadrature. The
pressure Poisson equation, and Helmholtz equation in case of the Stokes equations, are solved with
the Conjugate Gradient method and a hybrid AMG-DG preconditioner. This particular DG multigrid
algorithm is based on a correction in the conforming piecewise linear subspace where only the low
order components are explicitly assembled and the operations on the DG level are done
matrix-free. In case of the Navier-Stokes equations, the arising equations in the viscous substep
are solved with matrix-free Newton-GMRes. Postprocessing steps like the pressure update that involve
a mass matrix only, are solved with the matrix-free inverse mass matrix operator. As preconditioners
in the GMRes method or as DG smoothers within the multigrid algorithm, respectively, we employ
block Jacobi, Gauss-Seidel or symmetric Gauss-Seidel methods for use in domain
decomposition. In order to solve the diagonal blocks, the following approaches have been
implemented. In the first (named partially matrix-free), these diagonal blocks are factorized and
during preconditioning a forward/backward solve is used. In the second approach the diagonal
blocks are solved iteratively using matrix-free sum-factorization. Both of these variants have been
developed in \cite{BASTIAN2019417}. The third variant, however, implements the tensor product
preconditioners \cite{PAZNER2018344} which are based on Kronecker singular value decomposition
(KSVD) of the diagonal blocks.

\subsection{Potential flow}
\label{sec:potential-flow}

We consider the potential flow \cite{doi:10.1137/17M1112017} of the form \(v(x,y,t) =
t\tnabla\chi(x,y)\) with the harmonic potential \(\chi(x,y) = 5 x^4 y + y^5 - 10 x^2 y^3\) such that \(\Delta\chi(x,y)
= 0\). By construction, \(v\) is divergence-free and also harmonic. The velocity field solves the
instationary Stokes equations on \(\Omega=(0,1)^2\)
\begin{align*}
  \partial_t v - \nu\Delta v + \tnabla p &= 0 \\
  \nabla\cdot v &= 0
\end{align*}
for any viscosity \(\nu\) with pressure \(p = -\chi\). Global Dirichlet boundary conditions are imposed
by the analytical solution. This test case serves as a investigation on pressure robustness for the
three variants of Helmholtz decomposition within the RIPCS. It is intended to demonstrate how the
velocity errors may depend (or may not depend) on the viscosity. For a pressure robust method, it is
expected that these errors are independent on the viscosity.

We simulate the potential flow problem starting from time \(t_0 = 0\) up to \(T = 1\) for different
viscosity values \(\nu\in\{10^{-1}, 10^{-2}, 10^{-3}, 10^{-4}, 10^{-5}\}\). We consider the cumulative
velocity errors \(L^2(t_0,T;L^2(\Omega)^d)\) and \(L^2(t_0,T;H_0^1(\Omega)^d)\), and the cumulative pressure error
\(L^2(t_0,T;L^2(\Omega))\) for varying spatial resolutions where the time integral is approximated by a
trapezoidal rule. Therefore we carry out computations on rectangular meshes with \(2^{l+2}\) cells
per direction where the level \(l\) ranges over \(\{1,...,5\}\). To eliminate the temporal error we
use a fixed time step size of \(\Delta t = 5\cdot 10^{-3}\). The upcoming numerical analysis is performed for
the polynomial degrees \(p=2,3\) as the corresponding velocity ansatz spaces do not represent the
exact solution yet.

\begin{figure}[!htb]
  \hspace*{-2em}
  \includegraphics[width=1.2\textwidth]{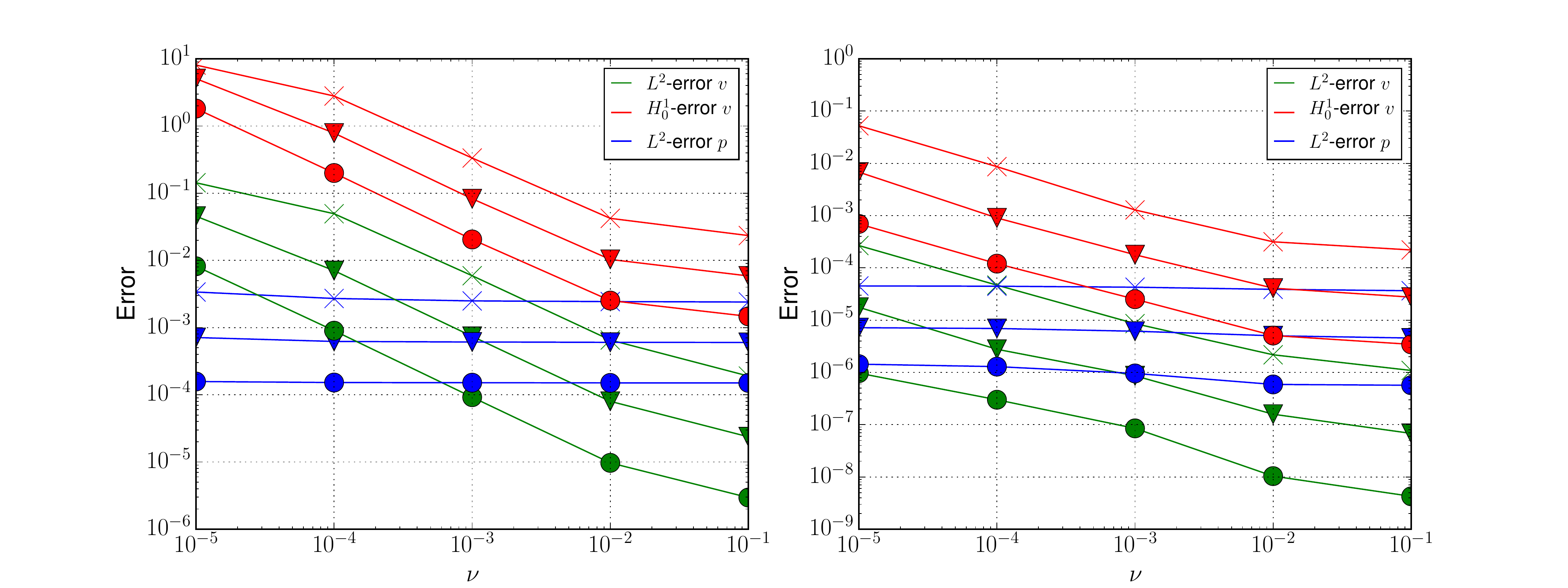}
  \caption{Errors accumulated over time vs. the viscosity parameter \(\nu\in\{10^{-5},...,10^{-1}\}\) for
    the div-div projection on different mesh refinement levels and fixed time step size \(\Delta t =
    5\cdot 10^{-3}\). \((\times)\) denotes the errors for 16 cells per direction (level 2) on the
    unitsquare,\((\blacktriangledown)\) for 32 cells per direction (level 3) and \((\bullet)\) for 64
    cells per direction (level 4).
    \newline
    Left subfigure shows results for \(p=2\), right subfigure \(p=3\).}
  \label{fig:invariance_property_divdiv_viscosity}
\end{figure}

\begin{figure}[!htb]
  \hspace*{-2em}
  \includegraphics[width=1.2\textwidth]{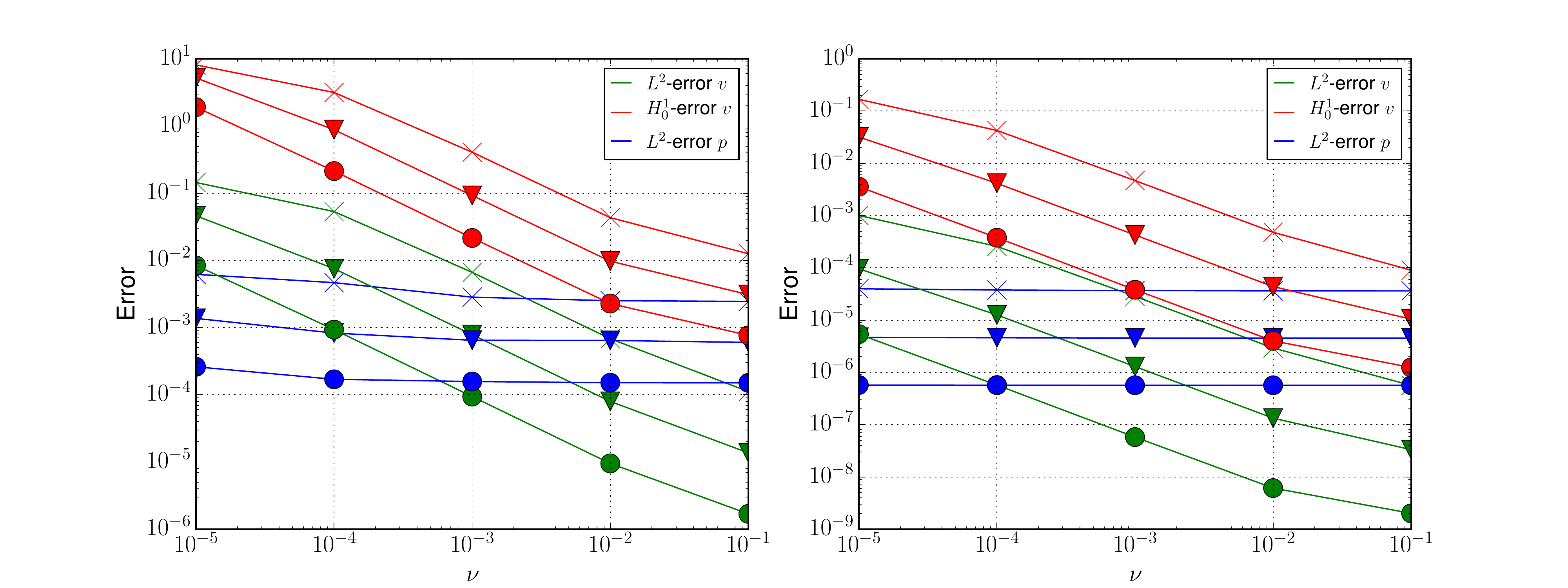}
  \caption{Errors accumulated over time vs. the viscosity parameter \(\nu\in\{10^{-5},...,10^{-1}\}\) for
    the pressure Poisson Raviart-Thomas projection in \(\RT_h^{p-2}\) on different mesh refinement
    levels and fixed time step size \(\Delta t = 5\cdot 10^{-3}\).
    \((\times)\) denotes the errors for 16 cells per direction (level 2) on the
    unitsquare,\((\blacktriangledown)\) for 32 cells per direction (level 3) and \((\bullet)\) for 64
    cells per direction (level 4).
    \newline
    Left subfigure shows results for \(p=2\), right subfigure \(p=3\).}
  \label{fig:invariance_property_pressurepoisson_viscosity}
\end{figure}

Table \ref{tab:invariance_property_divdiv_spatialrates} shows the convergence behavior in space for
\(\nu = 10^{-1}\) obtained by the div-div projection. The order of convergence in \(L^2\) and
\(H^1\) for the velocity, respectively, and for the pressure in \(L^2\) is optimal as predicted by
classical theory. The rates do not deteriorate for smaller viscosity values, though the velocity
errors exhibit a growth proportional to the factor \(1/\nu\) that comes with the bound including the
pressure as can be seen in figure \ref{fig:invariance_property_divdiv_viscosity}. In all the figures
presented the green curves show the velocity \(L^2\)-error for 16, 32, 64 cells per direction. The
curves in red color show the velocity \(H_0^1\)-error, the curves in blue color the pressure
\(L^2\)-error for the same spatial resolutions. The numerical results thus indicate that the div-div
projection utilized within the RIPCS does not provide a pressure robust splitting scheme. The same
conclusion can be drawn from pressure Poisson Raviart-Thomas reconstruction, c.f. with the Table
\ref{tab:invariance_property_pressurepoisson_spatialrates} and figure
\ref{fig:invariance_property_pressurepoisson_viscosity}.
Although not explicitly shown here, we have performed the numerical experiments with the
div-div-conti projection. In \cite{2018CMAME.341..917A} it was proven for mixed DG discretizations
that the viscosity growth of the pressure term in the velocity estimate is reduced to
\(\nu^{-\frac12}\), if both stabilization terms are simultaneously added on quadrilateral/hexahedral
meshes. In fact, we observe an increase of the velocity errors close to \(\nu^{-\frac12}\) for
decreasing \(\nu\).
%
%
\begin{table}[H]
  \centering
  \begin{filecontents*}{pictures/invariance-property/invariance-property-divdiv/spatial_rates_mu1_sorder2.csv}
level, $L^2$-error $v$, $L^2$-rate $v$, $H_0^1$-error $v$, $H_0^1$-rate $v$, $L^2$-error $p$, $L^2$-rate $p$ 
1 , 1.518e-03,  , 9.307e-02,  , 9.681e-03, 
2 , 1.901e-04, 2.997, 2.343e-02, 1.990, 2.406e-03, 2.009
3 , 2.371e-05, 3.003, 5.871e-03, 1.996, 6.013e-04, 2.000
4 , 2.959e-06, 3.002, 1.469e-03, 1.999, 1.503e-04, 2.000
5 , 3.695e-07, 3.001, 3.675e-04, 1.999, 3.758e-05, 2.000
\end{filecontents*}
\pgfplotstabletypeset[
  col sep=comma,
  string type,
  display columns/1/.style={column type={c}},
  display columns/2/.style={column type={|c}},
  display columns/3/.style={column type={|c}},
  display columns/4/.style={column type={|c}},
  display columns/5/.style={column type={|c}},
  display columns/6/.style={column type={|c}},
  every head row/.style={after row=\hline},
  every last row/.style={after row={\hline\medskip}},
  ]{pictures/invariance-property/invariance-property-divdiv/spatial_rates_mu1_sorder2.csv}
  \begin{filecontents*}{pictures/invariance-property/invariance-property-divdiv/spatial_rates_mu1_sorder3.csv}
level, $L^2$-error $v$, $L^2$-rate $v$, $H_0^1$-error $v$, $H_0^1$-rate $v$, $L^2$-error $p$, $L^2$-rate $p$ 
1 , 1.723e-05,  , 1.740e-03,  , 2.927e-04, 
2 , 1.089e-06, 3.984, 2.201e-04, 2.983, 3.653e-05, 3.002
3 , 6.807e-08, 4.000, 2.752e-05, 3.000, 4.544e-06, 3.007
4 , 4.254e-09, 4.000, 3.439e-06, 3.000, 5.674e-07, 3.001
5 , 2.661e-10, 3.999, 4.299e-07, 3.000, 7.091e-08, 3.000
\end{filecontents*}
\pgfplotstabletypeset[
  col sep=comma,
  string type,
  display columns/1/.style={column type={c}},
  display columns/2/.style={column type={|c}},
  display columns/3/.style={column type={|c}},
  display columns/4/.style={column type={|c}},
  display columns/5/.style={column type={|c}},
  display columns/6/.style={column type={|c}},
  every head row/.style={after row=\hline},
  every last row/.style={after row=\hline},
  ]{pictures/invariance-property/invariance-property-divdiv/spatial_rates_mu1_sorder3.csv}
  \caption{Errors accumulated over time and convergence rates for the potential flow problem with
    \(\nu = 10^{-1}\) utilizing div-div projection. The top shows results for polynomial degree 2, the
    bottom for polynomial degree 3.}
  \label{tab:invariance_property_divdiv_spatialrates}
\end{table}

%
%
\begin{table}[H]
  \centering
  \begin{filecontents*}{pictures/invariance-property/invariance-property-pressurepoisson-rt/spatial_rates_mu1_sorder2.csv}
level, $L^2$-error $v$, $L^2$-rate $v$, $H_0^1$-error $v$, $H_0^1$-rate $v$, $L^2$-error $p$, $L^2$-rate $p$ 
1 , 9.165e-04,  , 5.165e-02,  , 9.874e-03, 
2 , 1.116e-04, 3.038, 1.256e-02, 2.039, 2.453e-03, 2.009
3 , 1.368e-05, 3.028, 3.092e-03, 2.023, 6.015e-04, 2.028
4 , 1.691e-06, 3.016, 7.667e-04, 2.012, 1.503e-04, 2.000
5 , 2.102e-07, 3.008, 1.909e-04, 2.006, 3.758e-05, 2.000
\end{filecontents*}
\pgfplotstabletypeset[
  col sep=comma,
  string type,
  display columns/1/.style={column type={c}},
  display columns/2/.style={column type={|c}},
  display columns/3/.style={column type={|c}},
  display columns/4/.style={column type={|c}},
  display columns/5/.style={column type={|c}},
  display columns/6/.style={column type={|c}},
  every head row/.style={after row=\hline},
  every last row/.style={after row={\hline\medskip}},
  ]{pictures/invariance-property/invariance-property-pressurepoisson-rt/spatial_rates_mu1_sorder2.csv}
  \begin{filecontents*}{pictures/invariance-property/invariance-property-pressurepoisson-rt/spatial_rates_mu1_sorder3.csv}
level, $L^2$-error $v$, $L^2$-rate $v$, $H_0^1$-error $v$, $H_0^1$-rate $v$, $L^2$-error $p$, $L^2$-rate $p$ 
1 , 9.862e-06,  , 8.080e-04,  , 2.902e-04, 
2 , 5.635e-07, 4.129, 9.016e-05, 3.164, 3.631e-05, 2.999
3 , 3.324e-08, 4.083, 1.042e-05, 3.113, 4.539e-06, 3.000
4 , 2.010e-09, 4.048, 1.242e-06, 3.068, 5.673e-07, 3.000
5 , 1.237e-10, 4.023, 1.513e-07, 3.037, 7.091e-08, 3.000
\end{filecontents*}
\pgfplotstabletypeset[
  col sep=comma,
  string type,
  display columns/1/.style={column type={c}},
  display columns/2/.style={column type={|c}},
  display columns/3/.style={column type={|c}},
  display columns/4/.style={column type={|c}},
  display columns/5/.style={column type={|c}},
  display columns/6/.style={column type={|c}},
  every head row/.style={after row=\hline},
  every last row/.style={after row=\hline},
  ]{pictures/invariance-property/invariance-property-pressurepoisson-rt/spatial_rates_mu1_sorder3.csv}
  \caption{Errors accumulated over time and convergence rates for the potential flow problem with
    \(\nu = 10^{-1}\) utilizing pressure Poisson Raviart-Thomas projection in \(\RT_h^{p-2}\). The
    top shows results for polynomial degree 2, the bottom for polynomial degree 3.}
  \label{tab:invariance_property_pressurepoisson_spatialrates}
\end{table}

%
%
\begin{table}[H]
  \centering
  \begin{filecontents*}{pictures/invariance-property/invariance-property-helmholtz-rt/spatial_rates_mu1_sorder2.csv}
level, $L^2$-error $v$, $L^2$-rate $v$, $H_0^1$-error $v$, $H_0^1$-rate $v$, $L^2$-error $p$, $L^2$-rate $p$ 
1 , 1.720e-02,  , 9.861e-01,  , 1.917e-02, 
2 , 4.014e-03, 2.100, 4.929e-01, 1.000, 2.653e-03, 2.853
3 , 9.950e-04, 2.012, 2.465e-01, 1.000, 6.085e-04, 2.124
4 , 2.486e-04, 2.001, 1.233e-01, 1.000, 1.508e-04, 2.013
5 , 6.216e-05, 2.000, 6.163e-02, 1.000, 3.760e-05, 2.004
\end{filecontents*}
\pgfplotstabletypeset[
  col sep=comma,
  string type,
  display columns/1/.style={column type={c}},
  display columns/2/.style={column type={|c}},
  display columns/3/.style={column type={|c}},
  display columns/4/.style={column type={|c}},
  display columns/5/.style={column type={|c}},
  display columns/6/.style={column type={|c}},
  every head row/.style={after row=\hline},
  every last row/.style={after row={\hline\medskip}},
  ]{pictures/invariance-property/invariance-property-helmholtz-rt/spatial_rates_mu1_sorder2.csv}
  \begin{filecontents*}{pictures/invariance-property/invariance-property-helmholtz-rt/spatial_rates_mu1_sorder3.csv}
level, $L^2$-error $v$, $L^2$-rate $v$, $H_0^1$-error $v$, $H_0^1$-rate $v$, $L^2$-error $p$, $L^2$-rate $p$ 
1 , 3.530e-04,  , 3.543e-02,  , 2.958e-04, 
2 , 4.357e-05, 3.018, 9.079e-03, 1.964, 3.761e-05, 2.975
3 , 5.447e-06, 3.000, 2.270e-03, 2.000, 4.686e-06, 3.005
4 , 6.798e-07, 3.002, 5.652e-04, 2.006, 5.718e-07, 3.035
5 , 8.496e-08, 3.000, 1.410e-04, 2.002, 7.100e-08, 3.010
\end{filecontents*}
\pgfplotstabletypeset[
  col sep=comma,
  string type,
  display columns/1/.style={column type={c}},
  display columns/2/.style={column type={|c}},
  display columns/3/.style={column type={|c}},
  display columns/4/.style={column type={|c}},
  display columns/5/.style={column type={|c}},
  display columns/6/.style={column type={|c}},
  every head row/.style={after row=\hline},
  every last row/.style={after row=\hline},
  ]{pictures/invariance-property/invariance-property-helmholtz-rt/spatial_rates_mu1_sorder3.csv}
  \caption{Errors accumulated over time and convergence rates for the potential flow problem with
    \(\nu = 10^{-1}\) utilizing Helmholtz flux Raviart-Thomas projection in \(\RT_h^{p-1}\). The top
    shows results for polynomial degree 2, the bottom for polynomial degree 3.}
  \label{tab:invariance_property_helmholtz_spatialrates}
\end{table}

\begin{figure}[!htb]
  \hspace*{-2em}
  \includegraphics[width=1.2\textwidth]{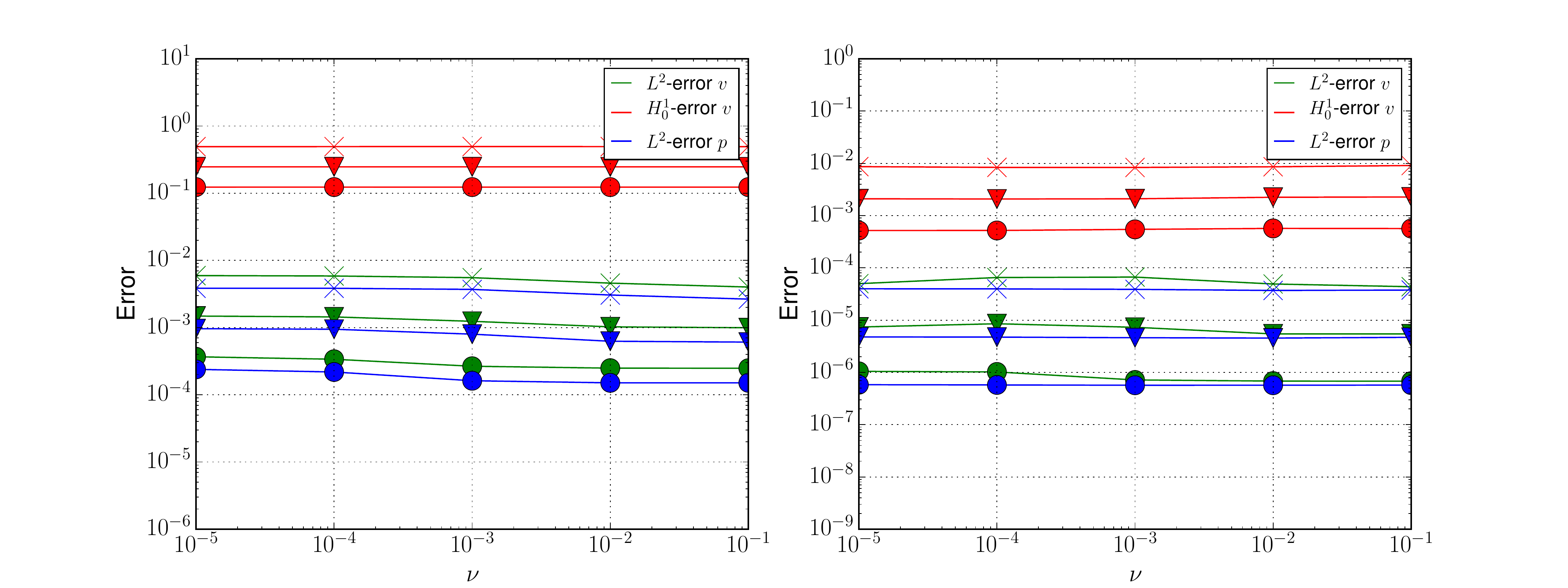}
  \caption{Errors accumulated over time vs. the viscosity parameter \(\nu\in\{10^{-5},...,10^{-1}\}\) for
    the Helmholtz flux Raviart-Thomas projection in \(\RT_h^{p-1}\) on different mesh refinement
    levels and fixed time step size \(\Delta t = 5\cdot 10^{-3}\).
    \((\times)\) denotes the errors for 16 cells per direction (level 2) on the
    unitsquare,\((\blacktriangledown)\) for 32 cells per direction (level 3) and \((\bullet)\) for 64
    cells per direction (level 4).
    \newline
    Left subfigure shows results \(p=2\), right subfigure \(p=3\).}
  \label{fig:invariance_property_helmholtz_viscosity}
\end{figure}

Figure \ref{fig:invariance_property_helmholtz_viscosity} shows the dependence of the errors on the
viscosity obtained by Helmholtz flux reconstruction, and Table
\ref{tab:invariance_property_helmholtz_spatialrates} the convergence behavior in space for \(\nu =
10^{-1}\). Importantly, as can be seen in the figure, the dependence of the velocity errors on the
viscosity is visually absent for Helmholtz flux reconstruction. The implementation, however,
currently computes a tentative velocity in \(X_h^3\) in order to reconstruct a solenoidal velocity in
\(\RT_h^2\), a subject that we will discuss in Section \ref{sec:conclusion-outlook}. The space
\(X_h^3\) delivers one additional approximation order for the velocity compared to \(\RT_h^2\) and
therefore we will compare the accuracy to \(p=2\) methods. Note that the approximation order for
the pressure is unaffected by Helmholtz flux reconstruction. The right subfigure of
\ref{fig:invariance_property_helmholtz_viscosity} shows the dependence of errors on \(\nu\) for
Helmholtz flux reconstruction in \(\RT_h^2\), the left subfigure of
\ref{fig:invariance_property_divdiv_viscosity} the dependence for the spatial discretization with
\(p=2\) and div-div projection employed. It can be seen that the \(H^1\)- and \(L^2\)-errors for those
two methods are of equal magnitude for the largest viscosity value \(\nu = 10^{-1}\). But eventually
for decreasing \(\nu\), Helmholtz flux reconstruction in \(\RT_h^2\) outperforms div-div projection
with polynomial degree \(p=2\) as the \(H^1\)- and \(L^2\)-errors exhibit no growth. The same
conclusion can be drawn when comparing to pressure Poisson flux reconstruction.

We have also performed the same tests on a parallelepiped domain \(\Omega\). We observe that the velocity
errors still remain constant in \(\nu\) for domains constructed by affine transformations. The
corresponding results have been left out since they do not differ significantly compared to those in
figure \ref{fig:invariance_property_helmholtz_viscosity}.  The numerical experiments thereby
indicate that Helmholtz flux Raviart-Thomas projection utilized within the RIPCS provides a pressure
robust splitting scheme.

\subsubsection*{Potential flow with irrotational force translation}

To see how a discretization scheme fulfills the invariance property, consider the potential flow
problem with a irrotational source term \(f = \tnabla \psi\). As proposed in \cite{FVCA8Benchmark},
we choose \(\psi(x,y) = \exp(-10(1-x+2y))\). With this source term the instationary Stokes equations
\begin{align*}
  \partial_t v - \nu\Delta v + \tnabla p &= f \\
  \nabla\cdot v &= 0 \; .
\end{align*}
have the exact solution \(v = t\tnabla\chi\) and \(p = -\chi + \psi\). We have repeated the numerical
experiments to verify if the splitting schemes offer a invariance property or how far they deviate
from it depending on the viscosity. The numerical results demonstrate that the div-div and pressure
Poisson Raviart-Thomas projection do not yield a pressure robust splitting scheme because the
velocity errors are constantly higher as for the zero right-hand side. Moreover, the velocity errors
show the same growth in inverse proportion to the viscosity. In contrast, the velocity errors for
the Helmholtz flux Raviart-Thomas projection appear to be robust in pressure and viscosity.

\subsection{Gresho vortex}
\label{sec:gresho-vortex}

\begin{figure}[!htb]
  \centering
  \includegraphics[width=0.6\textwidth]{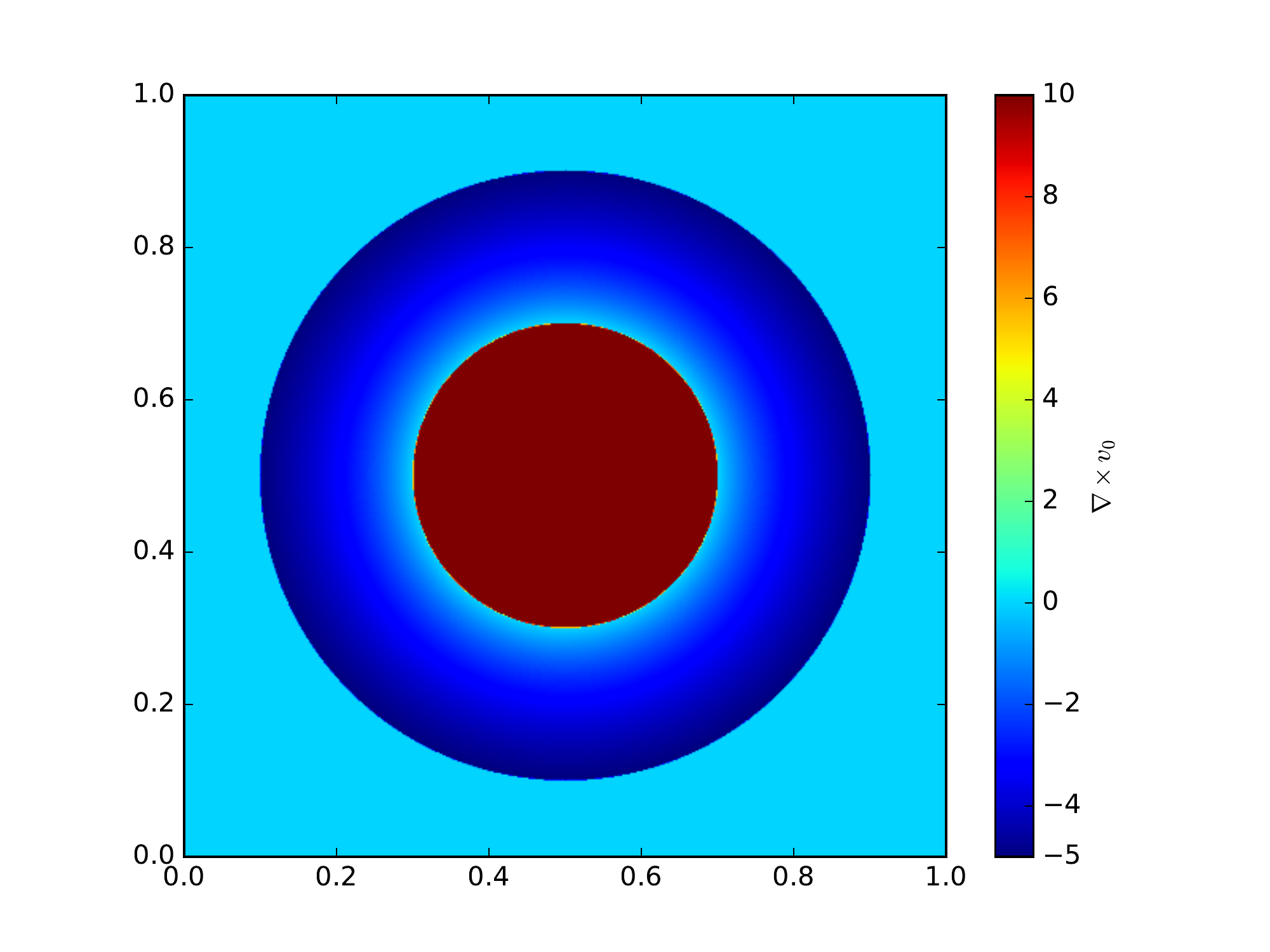}
  \caption{Vorticity of the initial Gresho vortex state.}
  \label{fig:gresho_initial_vorticity}
\end{figure}

The Gresho vortex has been recently proposed as model problem in \cite{GAUGER201989} for
investigating how well a discretization scheme preserves structures. It is argued that a pressure
robust method in this regard is in general superior to a non pressure robust method such that the
latter has certain difficulties preserving the vortex structure of this flow. Centered at
\(c = (0.5, 0.5)^T\) with constant translational velocity \(w_0\in\mathbb{R}^2\) (called wind), the
setup is described by the initial condition
\begin{equation}
  v_0(x) = w_0 +
  \begin{cases}
    (-5\tilde x_2, 5\tilde x_1)^T &, 0\leq r < 0.2 \\
    \left(-\frac{2\tilde x_2}{r} + 5\tilde x_2, \frac{2\tilde x_1}{r} - 5\tilde x_1\right)^T &, 0.2\leq r < 0.4 \\
    0 &, 0.4\leq r
  \end{cases}
  \label{eq:gresho_initial_condition}
\end{equation}
with shifted coordinates \(\tilde x = x - c\) towards the center, \(r = \left\|\tilde x\right\|_2\).
The initial Gresho vortex satisfies \(\nabla\cdot v_0 = 0\). In the core \(r<0.2\) the vorticity \(\nabla\times v_0\)
is constantly 10.0, decreases linearly in the circular ring \(0.2\leq r<0.4\), and vanishes for \(r\geq 0.4\).
The initial vorticity is displayed in figure \ref{fig:gresho_initial_vorticity} where the
discontinuities at the two aforementioned shells are visible.
The initial state is evolved in time by the Navier-Stokes equations
\begin{align*}
  \partial_t v - \nu\Delta v + (v\cdot\nabla)v + \tnabla p &= 0 \\
  \nabla\cdot v &= 0
\end{align*}
in the periodic square \(\Omega = (0,1)^2\). As in \cite{GAUGER201989} we set the viscosity to
\(\nu = 10^{-5}\). For the \emph{standing} Gresho vortex with \(w_0 = 0\) the convection term
\((v_0\cdot\nabla)v_0\) can be balanced by the gradient of a pressure \(p_0\), \((v_0\cdot\nabla)v_0 = -\tnabla p_0\),
and, consequently, the standing Gresho vortex describes a steady solution to the instationary Euler
equations. With the diffusion term \(-\nu\Delta v\) added, the Navier-Stokes problem smoothes out the
discontinuities from the initial condition. Now, the choice \(w_0 = (1/3,1/3)^T\) is referred to as
the \emph{moving} Gresho vortex that we will consider in this work. In this case the vorticity
distribution is transported in the top-right direction through the periodic domain such that the
vortex at time \(T=3\) is intended to be again centered around \(c=(0.5,0.5)^T\).

We simulate the moving Gresho vortex problem from time \(t_0 = 0\) up to \(T = 3\) utilizing the
three variants of discrete Helmholtz decomposition within the RIPCS for the polynomial degrees 4 and
8. For \(p=4\) we use a spatial resolution of 32 cells per direction, for \(p=8\) 16 cells per
direction. All runs have been performed with a fixed time step size of \(\Delta t = 2\cdot 10^{-3}\).

\begin{figure}[!htb]
  \hspace*{-2em}
  \includegraphics[width=1.2\textwidth]{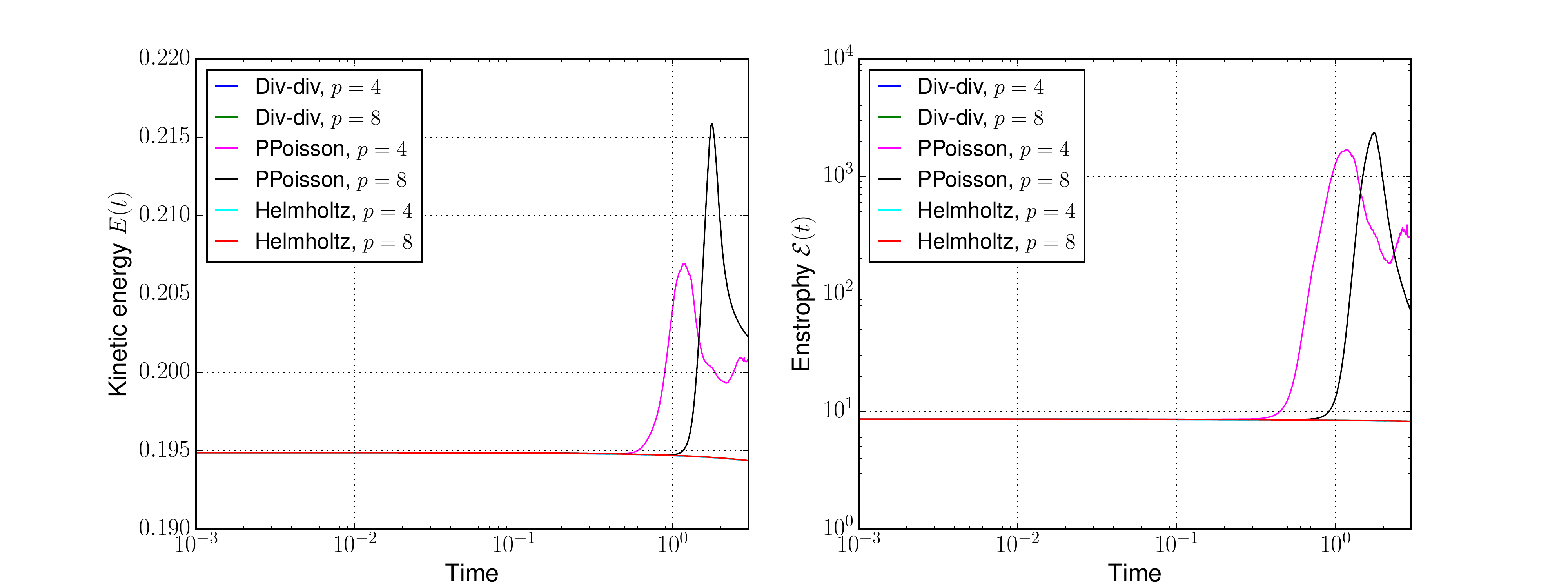}
  \caption{Evolution of kinetic energy (left) and enstrophy (right) for the \THcube{p}{p-1}
    pairs \((p\in\{4,8\})\) and respective discrete Helmholtz decompositions. \texttt{Helmholtz} in the
    legend abbreviates the reconstruction of the Helmholtz flux in \(\RT_h^{p-1}\), and
    \texttt{PPoisson} the reconstruction of the pressure Poisson flux in \(\RT_h^{p-2}\). Note that
    the evolution obtained by div-div projection and Helmholtz flux reconstruction are almost
    indistinguishable.}
  \label{fig:gresho_ekin_enstrophy}
\end{figure}

Let us start the discussion on the numerical results by considering important flow quantities for
each of the simulation performed. In figure \ref{fig:gresho_ekin_enstrophy} the kinetic energy \(E(t)
= \frac12 (v_h,v_h)_{0,\Omega}\) and enstrophy \(\mathcal{E}(t) = \frac12 (\nabla_h\times v_h,\nabla_h\times v_h)_{0,\Omega}\) are
displayed over time.
At the beginning all kinetic energy curves remain very close to constant as expected for this low
viscosity value. However, towards the end, the computations performed with the pressure Poisson
Raviart-Thomas reconstruction show an increase which is nonphysical for this freely evolving
system. The same observation can be made for the enstrophy which in theory cannot increase as well
in two dimensions. The div-div projection and Helmholtz flux Raviart-Thomas reconstruction deliver
stable results which are close to the Scott-Vogelius computations in \cite{GAUGER201989}.
This comparison shows that in terms of conserved quantities for inviscid flows, the splitting scheme
realized by pressure Poisson Raviart-Thomas projection is outperformed by the other two variants.

\begin{figure}[!htb]
  \includegraphics[width=\textwidth]{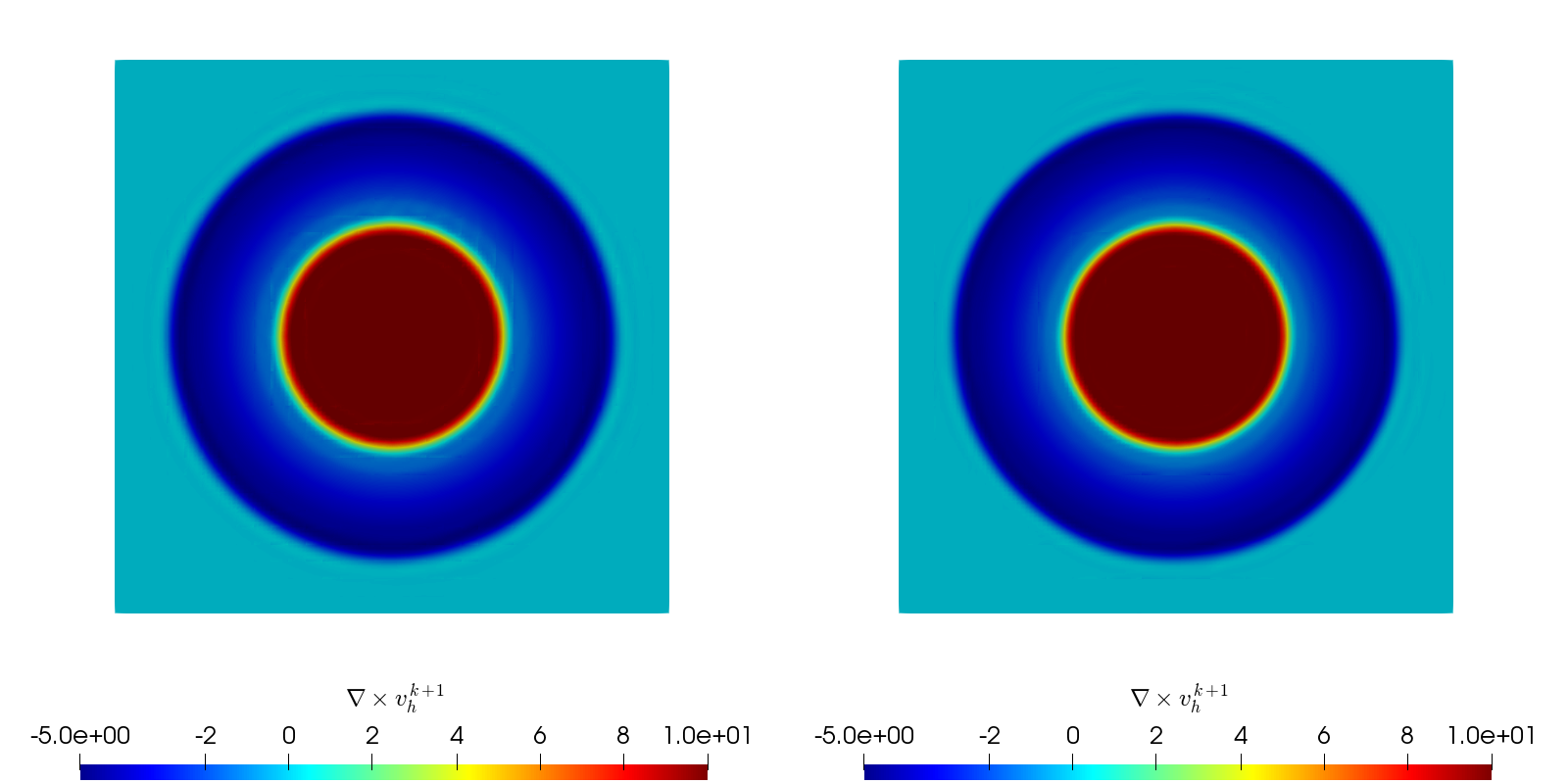}
  \caption{Vorticity of moving Gresho vortex at time \(T=3\) by \THcube{p}{p-1} discretization with
    div-div projection. Left part shows \(p=4\) method (32 cells per direction), right part the
    \(p=8\) method (16 cells per direction).}
  \label{fig:gresho_divdiv_final_vorticity}
\end{figure}

\begin{figure}[!htb]
  \includegraphics[width=\textwidth]{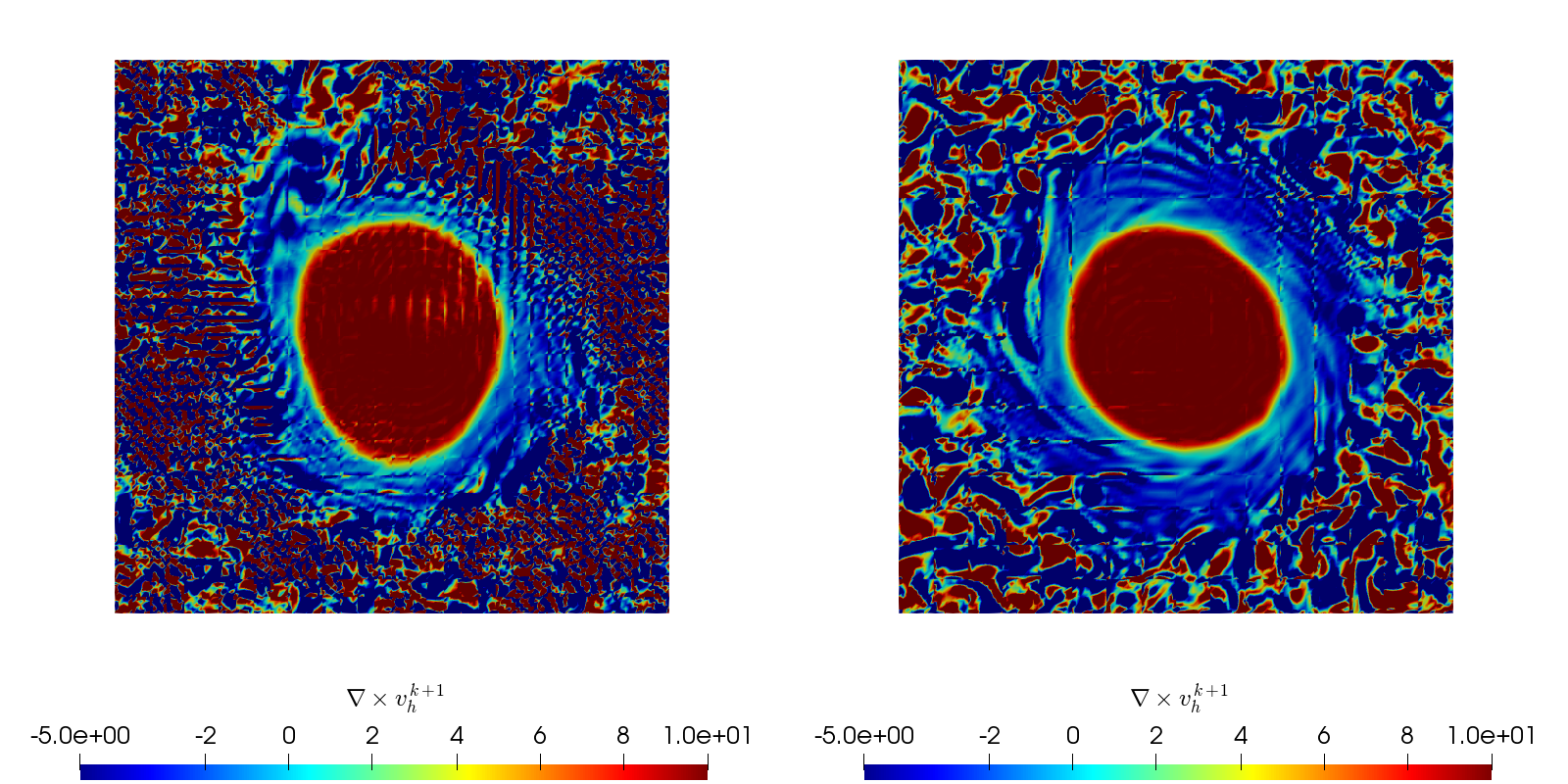}
  \caption{Vorticity of moving Gresho vortex at time \(T=3\) by \THcube{p}{p-1} discretization with
    pressure Poisson Raviart-Thomas reconstruction in \(\RT_h^{p-2}\). Left part shows \(p=4\)
    method (32 cells per direction), right part the \(p=8\) method (16 cells per direction).}
  \label{fig:gresho_pressurepoisson_final_vorticity}
\end{figure}

\begin{figure}[!htb]
  \includegraphics[width=\textwidth]{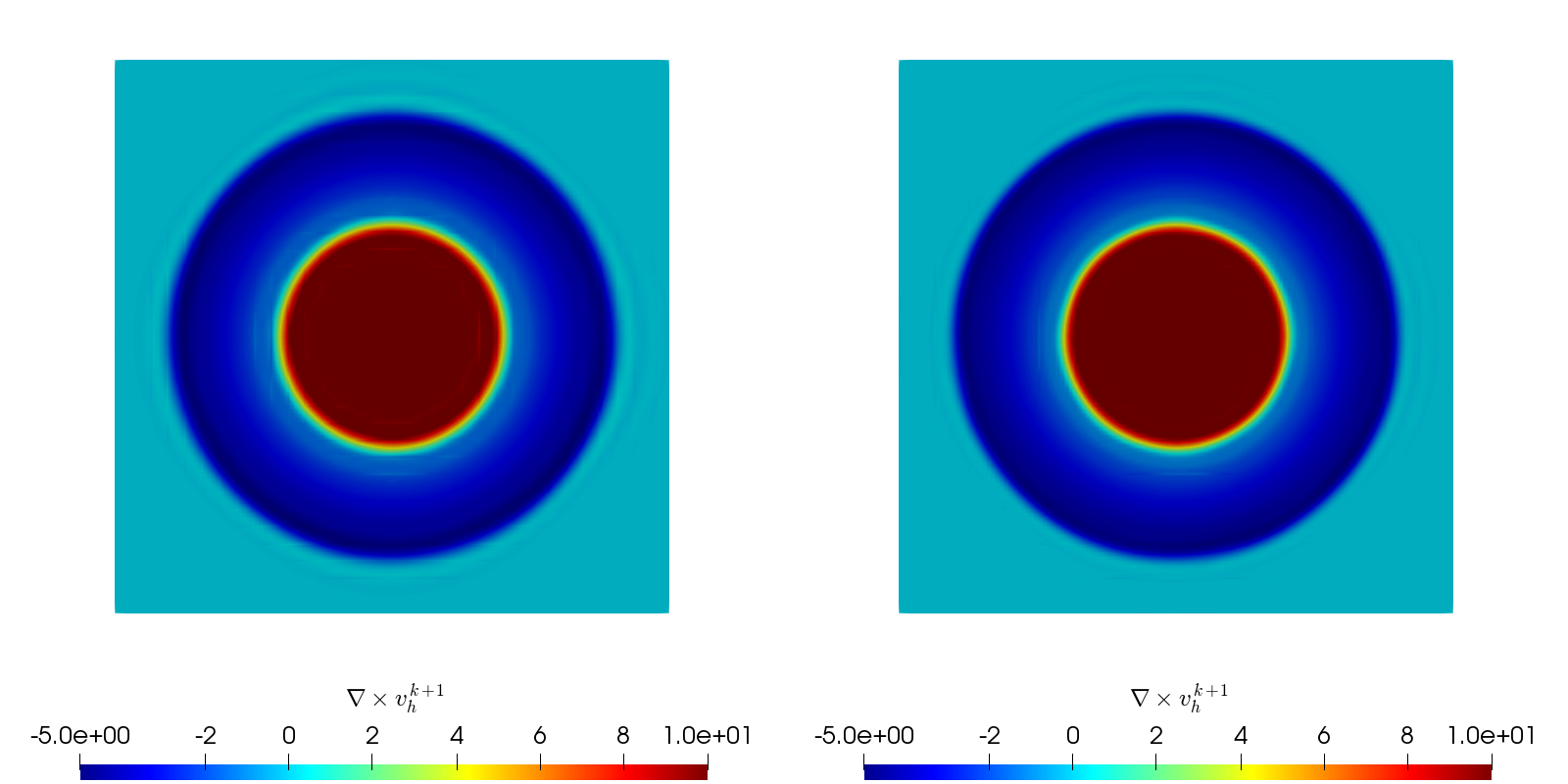}
  \caption{Vorticity of moving Gresho vortex at time \(T=3\) by \THcube{p}{p-1} discretization with
    Helmholtz flux Raviart-Thomas projection in \(\RT_h^{p-1}\). Left part shows \(p=4\) method (32
    cells per direction), right part the \(p=8\) method (16 cells per direction).}
  \label{fig:gresho_helmholtz_final_vorticity}
\end{figure}

Now, let us investigate the vorticity plots at time \(T=3\) for the three discrete Helmholtz
decompositions. Figure \ref{fig:gresho_divdiv_final_vorticity}-~\ref{fig:gresho_helmholtz_final_vorticity}
show the results for the polynomial degree 4 (left subfigures) and \(p=8\) (right subfigures). One
can observe that div-div projection and Helmholtz flux Raviart-Thomas projection are able to
preserve the structure of the initial condition whereas the simulations with pressure Poisson
reconstruction are interspersed with instabilities away from the vortex core. These numerical
artifacts are eventually responsible for the nonphysical growth in kinetic energy and
enstrophy. Furthermore, one can see that the higher order method with \(p=8\) gives for both stable
variants slightly better results in terms of over-/under-shoots. On the one hand, this is surprising
as the initial condition is discontinuous and the flow is governed by inviscid vortex transport, but
on the other hand viscous diffusion - no matter how small - eventually provides sufficient
regularity of the problem. As above we can conclude that the non pressure robust splitting scheme
with pressure Poisson reconstruction is inferior in preserving structures with large gradient
parts. We also want to point out that the effect of div-div stabilization gives similar results to
those obtained by a pressure robust method, c.f. figures \ref{fig:gresho_divdiv_final_vorticity} and
\ref{fig:gresho_helmholtz_final_vorticity}. This observation has been made with respect to grad-div
stabilization in the literature.

\subsection{Beltrami flow}
\label{sec:beltrami-flow}

The Beltrami flow is one of the rare test problems that describes a fully three-dimensional solution
of the instationary Navier-Stokes equations. It has been derived as a class of analytical solutions
in \cite{EthierSteinmann1994} by separation of variables. The velocity field and pressure read
\begin{equation}
\begin{split}
  v_1(x,y,z,t) &= -a\,e^ {-\nu d^2\,t }\,\left(e^{a\,x}\,\sin \left(d\,z+a\,y\right)+
 \cos \left(d\,y+a\,x\right)\,e^{a\,z}\right) \\
  v_2(x,y,z,t) &= -a\,e^ {-\nu d^2\,t }\,\left(e^{a\,x}\,\cos \left(d\,z+a\,y\right)+e^{
 a\,y}\,\sin \left(a\,z+d\,x\right)\right) \\
  v_3(x,y,z,t) &= -a\,e^ {-\nu d^2\,t }\,\left(e^{a\,y}\,\cos \left(a\,z+d\,x\right)+
 \sin \left(d\,y+a\,x\right)\,e^{a\,z}\right) \\
  p(x,y,z,t) &= p_0(t) - 0.5\,a^2\,\rho\,e^ {-2\nu d^2\,t }\,(2\,\cos \left(d\,y+a\,x
 \right)\,e^{a\,\left(z+x\right)}\,\sin \left(d\,z+a\,y\right) \\
  &\quad +2\,e^{ a\,\left(y+x\right)}\,\sin \left(a\,z+d\,x\right)\,\cos \left(d\,z+a
 \,y\right) \\
  &\quad +2\,\sin \left(d\,y+a\,x\right)\,e^{a\,\left(z+y\right)}\,
 \cos \left(a\,z+d\,x\right)+e^{2\,a\,z}+e^{2\,a\,y}+e^{2\,a\,x} )
\end{split}
\label{eq:beltrami_3dproblem}
\end{equation}
where $p_0(t)$ is chosen to ensure $\int_\Omega p(x,y,z,t) \mathrm{d} x \mathrm{d} y \mathrm{d} z = 0$ over time.
The Beltrami flow solves the instationary Navier-Stokes equations for any positive viscosity $\nu > 0$.
The computational domain is given by $\Omega = (-1,1)^3$. On \(\partial\Omega\) Dirichlet boundary conditions for the
velocity are imposed by the exact solution. The Beltrami flow has the property that the velocity and
vorticity vectors are aligned, i.e. \(d\; v - \nabla\times v = 0\), and that the nonlinear convection term is
balanced by the pressure gradient, \(\rho(v\cdot\nabla)v = -\tnabla p\). It is therefore interesting to see how
a discretization handles the irrotational part for \(\nu\ll 1\).

\begin{figure}[!htb]
  \hspace*{-2em}
  \includegraphics[width=1.2\textwidth]{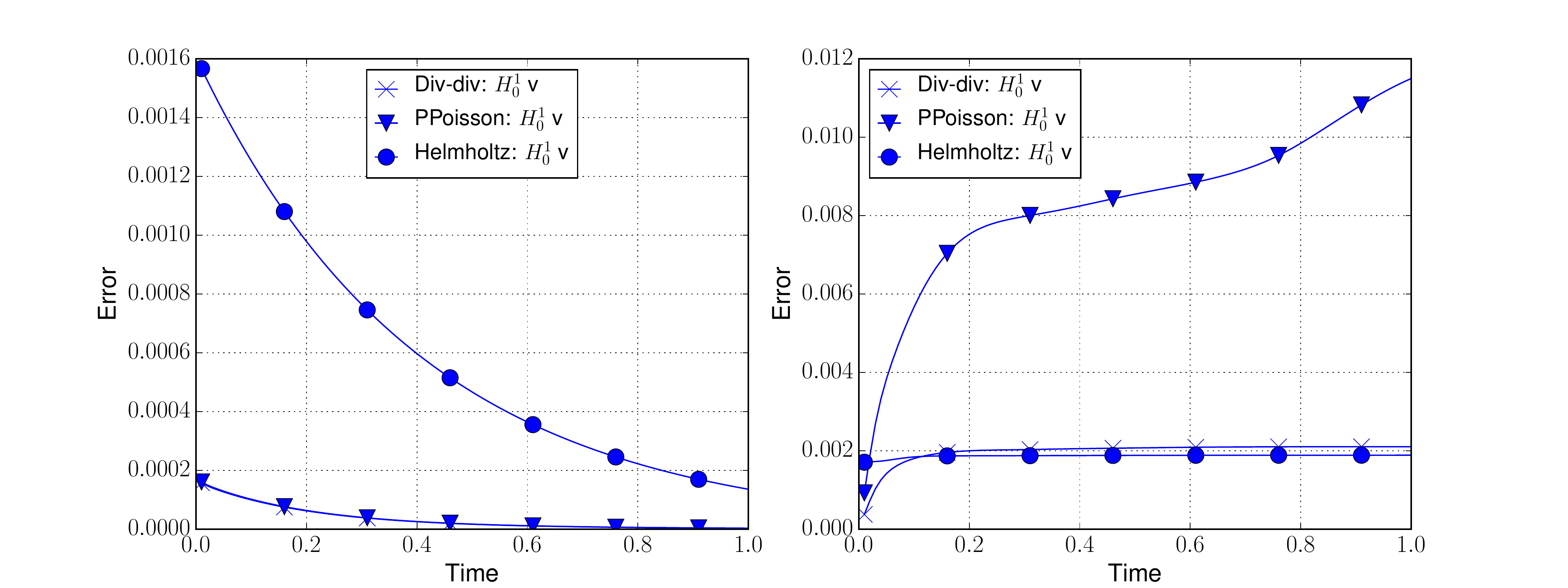}
  \caption{Plots of time versus \(H_0^1(\Omega)\)-velocity error for the three variants of discrete
    Helmholtz decomposition. Left part shows progression for \(\nu = 1\), right part for
    \(\nu = 10^{-3}\)}
  \label{fig:beltrami3d_spatialerrors_overtime}
\end{figure}

Inspired by the tests in \cite{LINKE2019350,GAUGER201989} we set the parameters
\(a=\pi/4,\;d=\pi/2\) and simulate the Beltrami problem for \(\nu\in\{1,10^{-3}\}\) from time \(t_0=0\) up to
\(T=1\). We set the polynomial degree for the velocity DG space to \(p=4\) and use a mesh
resolution of 8 cells per direction. To eliminate instabilities from the time-stepping we set
\(\Delta t = 10^{-3}\). Figure \ref{fig:beltrami3d_spatialerrors_overtime} shows the
\(|(v-v_h)(t)|_{1,\Omega}\) errors over time for the three variants of discrete Helmholtz
decomposition. The progression of the \(L^2\)-norm only differs in absolute values and has been
omitted.
One can observe that with \(\nu=1\) all splitting methods remain accurate (at the level of spatial
approximation error) up to the end time.
The right subfigure shows the corresponding results for viscosity \(10^{-3}\) where no such rapid
exponential decay of the errors is expected. Note that the initial errors are the same for all runs
but the y-scale has been enlarged accordingly. However, we observe that the \(H^1\)-seminorm errors
for pressure Poisson Raviart-Thomas projection grow in time and the solution becomes inaccurate. As
stated in the previous section, div-div projection improves the results such that only at the
beginning an increase is visible. In contrast, Helmholtz flux Raviart-Thomas projection within the
RIPCS provides constant \(|(v-v_h)(t)|_{1,\Omega}\) errors and, apart from the first time steps,
outperforms the div-div projection in this measure. We conclude that Helmholtz flux reconstruction
which is indicated to give a pressure robust splitting scheme reproduces well this (seemingly) easy
flow problem.

\subsection{3D Taylor-Green vortex}
\label{sec:3d-taylor-green}

The three-dimensional Taylor-Green vortex is aimed at testing the accuracy and performance of
high-order methods in a DNS. The flow starts from a simple large scale initial condition. In the
early phase it undergoes vortex stretching until laminar breakdown before the maximum dissipation
of the fluid is reached. The flow then transitions to turbulence followed by a decay phase of
eventually Homogeneous Isotropic Turbulence (HIT). The problem originates from
\cite{10.2307/96892} where classes of sinusoidal fields were considered as an initial condition that
satisfy the continuity equation \(\nabla\cdot v = 0\). It has been proposed as a reference benchmark since
the first edition of the international workshop on High-Order CFD Methods and is quoted in
\cite{TaylorGreen3D}, C3.5, for instance.

The simulation domain is \(\Omega=(-\pi L, \pi L)^3\) with periodic boundary conditions in all directions and
no external forcing, \(f=0\). The initial flow field is given by
\begin{align}
  \label{eq:taylorgreen_vortex3d_initial_v}
  \begin{split}
    v_0(x)\cdot e_1 &= V_0 \sin\left(\frac{x_1}{L}\right) \cos\left(\frac{x_2}{L}\right) \cos\left(\frac{x_3}{L}\right) \\
    v_0(x)\cdot e_2 &= -V_0 \cos\left(\frac{x_1}{L}\right) \sin\left(\frac{x_2}{L}\right) \cos\left(\frac{x_3}{L}\right) \\
    v_0(x)\cdot e_3 &= 0
  \end{split} \\
  \label{eq:taylorgreen_vortex3d_initial_p}
  p_0(x) &= p_0
  + \frac{\rho_0 V_0^2}{16} \left(\cos\left(\frac{2x_1}{L}\right) + \cos\left(\frac{2x_2}{L}\right)\right)
  \left(\cos\left(\frac{2x_3}{L}\right) + 2\right) \; .
\end{align}
The Reynolds number of the flow here is defined as \(\text{Re} = \frac{\rho_0 V_0 L}{\mu}\). As in the
references \cite{TaylorGreen3D,fenics:book} we set \(L=1, \; V_0 = 1, \; \rho_0 = \rho = 1, \; p_0 = 0, \;
\textup{Re} = 1600\).

\begin{figure}[!htb]
  \includegraphics[width=\textwidth]{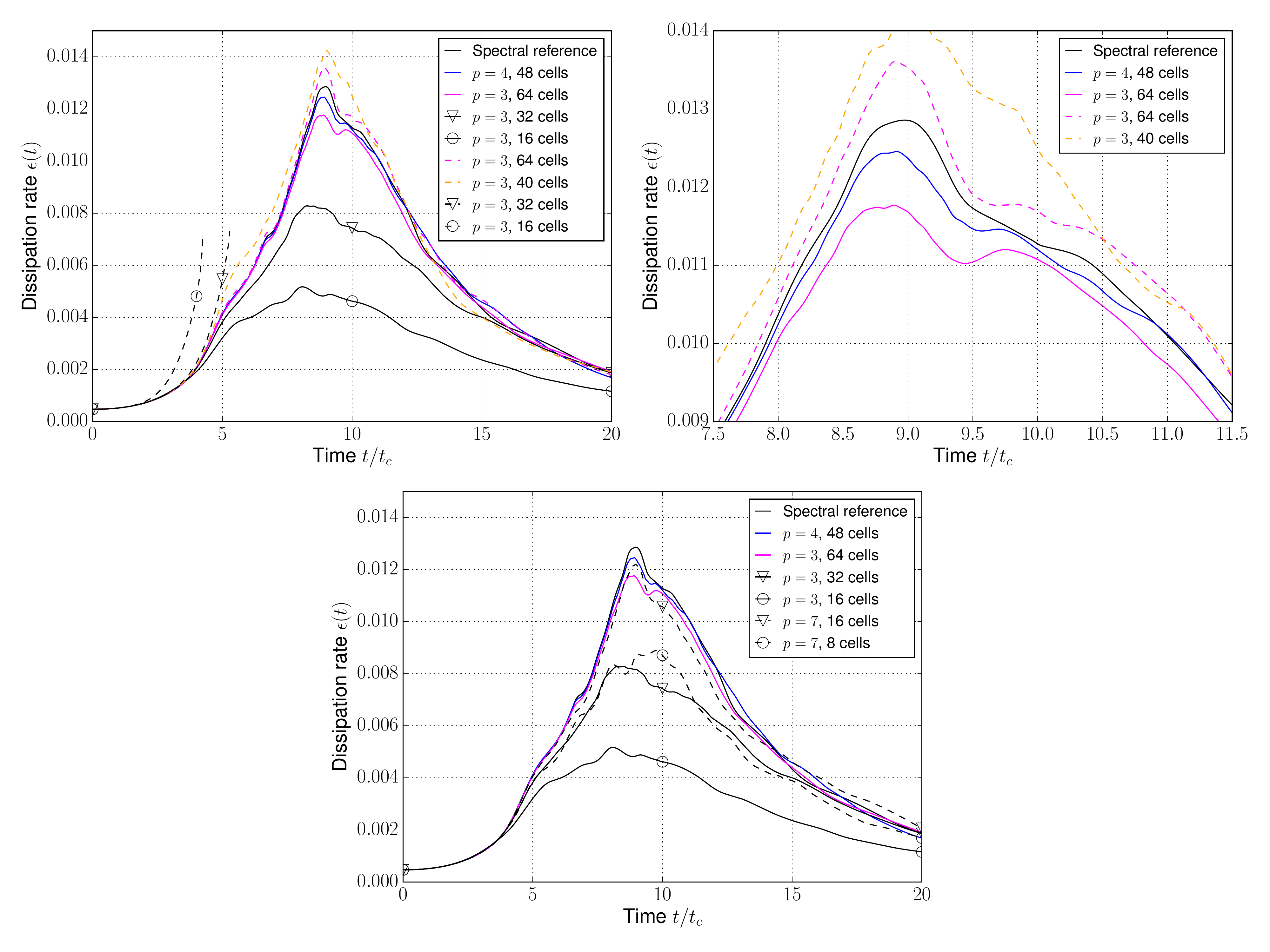}
  \caption{Time evolution of kinetic energy dissipation for the 3D Taylor-Green vortex. Comparison
    of various DG configurations with reference solution: In the upper plots the dashed lines refer
    to pressure Poisson Raviart-Thomas projection, and the solid lines using Helmholtz flux
    reconstruction. In the bottom center plot the solid lines correspond to the same configuration
    as above, but the dashed lines now correspond to a \(p=7\) run with Helmholtz flux
    reconstruction employed as well.}
  \label{fig:taylorgreen3d_diss_curves}
\end{figure}

Using this problem we want to demonstrate that with a pressure robust DG method, a robust method for
underresolved turbulent incompressible flows can be realized. To this end, we have done computations
on a series of uniformly refined, equidistant cuboid meshes for different polynomial degrees. There
are detailed reference data available which contain the temporal evolution of
\begin{itemize}
\item kinetic energy \(E(t) = \frac{1}{\rho_0 \left|\Omega\right|} \frac12 (\rho v_h,v_h)_{0,\Omega}\),
\item dissipation rate \(\epsilon(t) = \frac{\nu}{\left|\Omega\right|}(\nabla_h v_h, \nabla_h v_h)_{0,\Omega}\),
\item enstrophy \(\mathcal{E}(t) = \frac{1}{\rho_0 \left|\Omega\right|} \frac12 (\rho\nabla_h\times v_h,\nabla_h\times v_h)_{0,\Omega}\)
\end{itemize}
in the time interval \([0,20 t_c]\) where \(t_c = \frac{L}{V_0}\) is the convective time unit. The
reference values were obtained with a dealiased pseudo-spectral code run on \(512^3\) grid, time
integration was performed with a low-storage three-step Runge-Kutta method and a time step \(\Delta t = 10^{-3}\,t_c\).
We will concentrate on the results of kinetic energy dissipation \(\epsilon(t)\), and the spectral
distribution of kinetic energy \(E(k,t)\) s.t. \(E(t) = \int_0^\infty E(k,t) \dl k\) at \(t=9 t_c\), close
to the maximum amount of dissipation. Table \ref{tab:grid_configs_taylorgreen3d} summarizes
configurations, including variant of discrete Helmholtz decomposition, number of cells per
direction, polynomial degree and velocity DOFs.
\begin{table}[H]
  \centering
  \begin{tabular}{cccc}
    Helmholtz flux reconstruction & Cells per direction & \(p\) & Velocity DOFs \\ \cline{2-4}
    & 8 & 7 & 64  \\
    & 16 & 3 & 64 \\
    & 16 & 7 & 128 \\
    & 32 & 3 & 128 \\
    & 48 & 4 & 240 \\
    & 64 & 3 & 256 \\
    & 128 & 3 & 512 \\ \cline{2-4}
    \\
    Pressure Poisson reconstruction & Cells per direction & \(p\) & Velocity DOFs \\ \cline{2-4}
    & 16 & 3 & 64 \\
    & 32 & 3 & 128 \\
    & 40 & 3 & 160 \\
    & 64 & 3 & 256\\
    & 128 & 3 & 512 \\ \cline{2-4}
  \end{tabular}
  \caption{Grid configurations for 3D Taylor-Green problem.}
  \label{tab:grid_configs_taylorgreen3d}
\end{table}
We compare the results with the reference solution in figure \ref{fig:taylorgreen3d_diss_curves}.
The top left shows dissipation curves for different grid sizes and polynomial degree 3, including a
fully resolved simulation with \(p=4\) and 240 velocity DOFs. The dashed lines here refer to
pressure Poisson Raviart-Thomas projection. It can be seen that the corresponding underresolved
simulations (dashed lines, black color) are unstable which is caused by a crash of the computation
during build-up phase. The runs with 40 and 64 cells per direction (dashed lines in orange and
magenta color) describe essentially resolved simulations and already capture the shape from the
reference solution. The progression of the curves are similar to the ones obtained in
\cite{FEHN2018667} by a standard DG discretization with no additional stabilization terms
(c.f. figure 6 in this publication). In contrast, Helmholtz flux Raviart-Thomas projection provides
successful underresolved simulations for the two configurations: \(p=3\), 16 and 32 cells per
direction. The solid lines show the progression using this variant and one can observe that kinetic
energy dissipation is stronger underpredicted, the lower the spatial resolution is (c.f. with figure
6 in \cite{FEHN2018667}). The lop left plot shows a zoom around the maximum amount of dissipation
where the underresolved runs have been left out.

The bottom center plot shows simulations with the Helmholtz flux reconstruction only employed. Here,
the dashed lines show the results of the order \(p=7\) method with the same number of velocity DOFs
as for the \(p=3\) underresolved computations. One can see that the \(p=7\) method with 64 velocity
DOFs gives results with the same accuracy as the \(p=3\) with 128 velocity DOFs. Moreover, the
\(p=7\) method with 128 velocity DOFs closely matches the simulations with \(p=4\), 240 velocity
DOFs, which itself is more accurate than the simulation with \(p=3\), 256 velocity DOFs.
Hence, we also observe the improved accuracy per degree of freedom for higher polynomial \(p\) which
motivates the use of high-order methods for smooth problems. We further point out that the behavior
of our dissipation curves for different \(p\) and grid size \(h\) matches the enstrophy curves
presented in \cite{PAZNER2018344}.

\begin{figure}[!htb]
  \hspace*{2em}
  \includegraphics[width=0.75\textwidth]{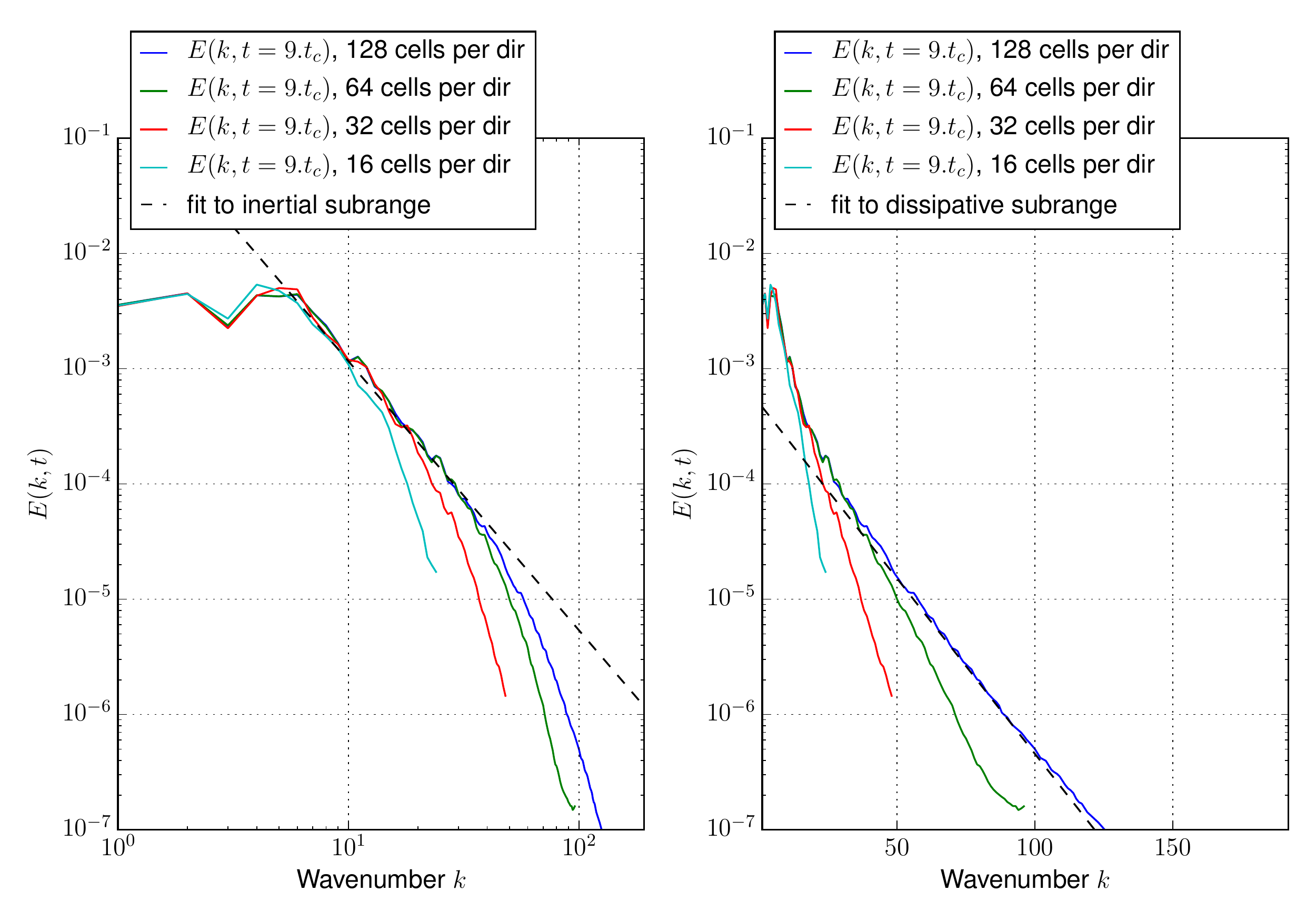}
  \caption{Kinetic energy spectra of the 3D Taylor-Green vortex at \(t=9t_c\) for polynomial degree
    \(p=3\) under \(h\)-coarsening, with fits to the inertial and dissipative subrange. All DG
    configurations here use Helmholtz flux reconstruction in \(\RT_h^{p-1}\).}
  \label{fig:taylorgreen3d_kinetic_spectra_hcoarsen}
\end{figure}

\begin{figure}[!htb]
  \hspace*{2em}
  \includegraphics[width=0.75\textwidth]{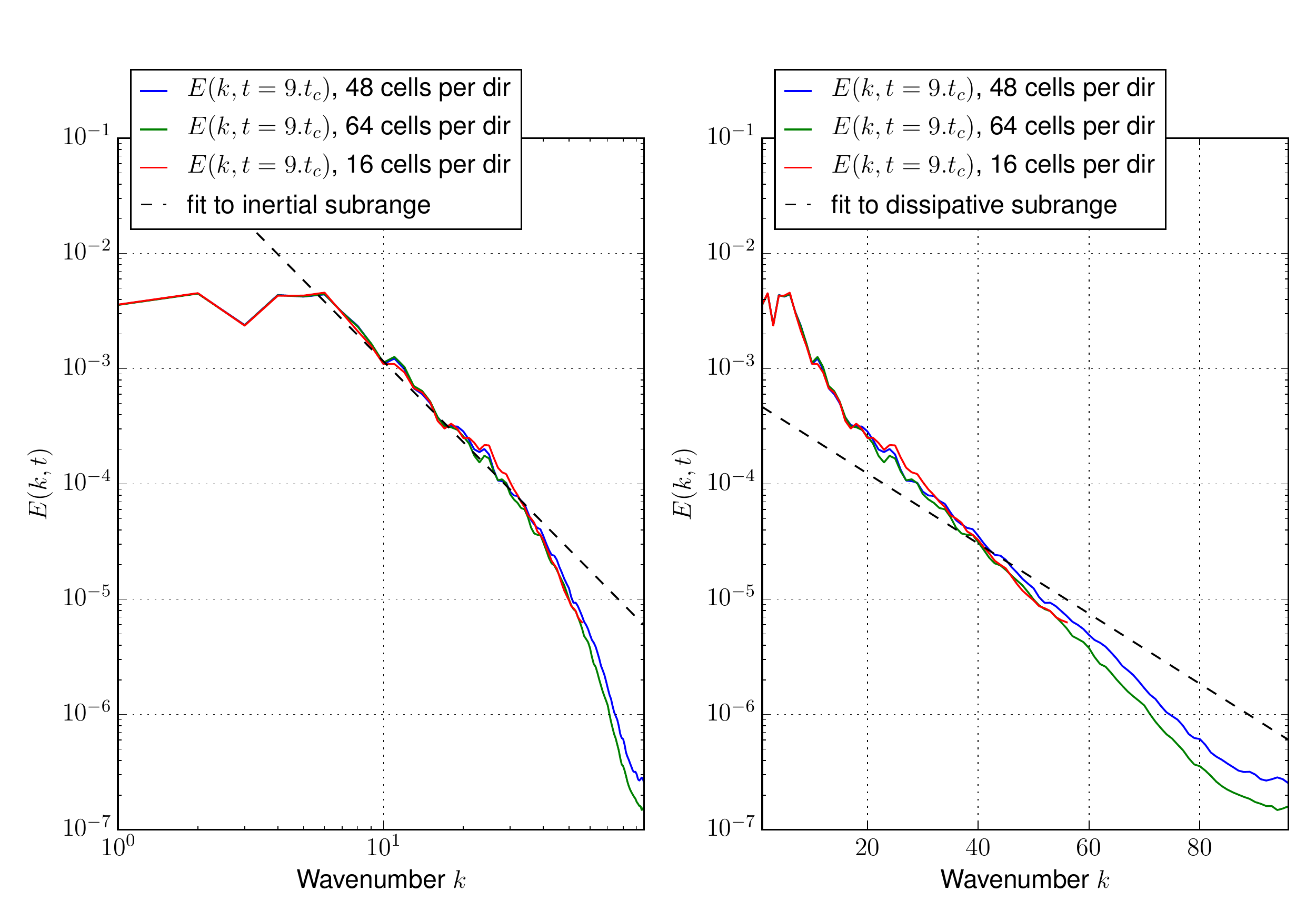}
  \caption{Kinetic energy spectra at \(t=9t_c\) under simultaneous \(p\)-refinement and
    \(h\)-coarsening. Green curve: \(p=3\) and 256 velocity DOFs. Blue curve: \(p=4\) and 240 velocity
    DOFs. Red curve: \(p=7\) and 128 velocity DOFs. All DG configurations here use Helmholtz flux
    reconstruction in \(\RT_h^{p-1}\).}
  \label{fig:taylorgreen3d_kinetic_spectra_prefine}
\end{figure}

\begin{figure}[!htb]
  \hspace*{2em}
  \includegraphics[width=0.75\textwidth]{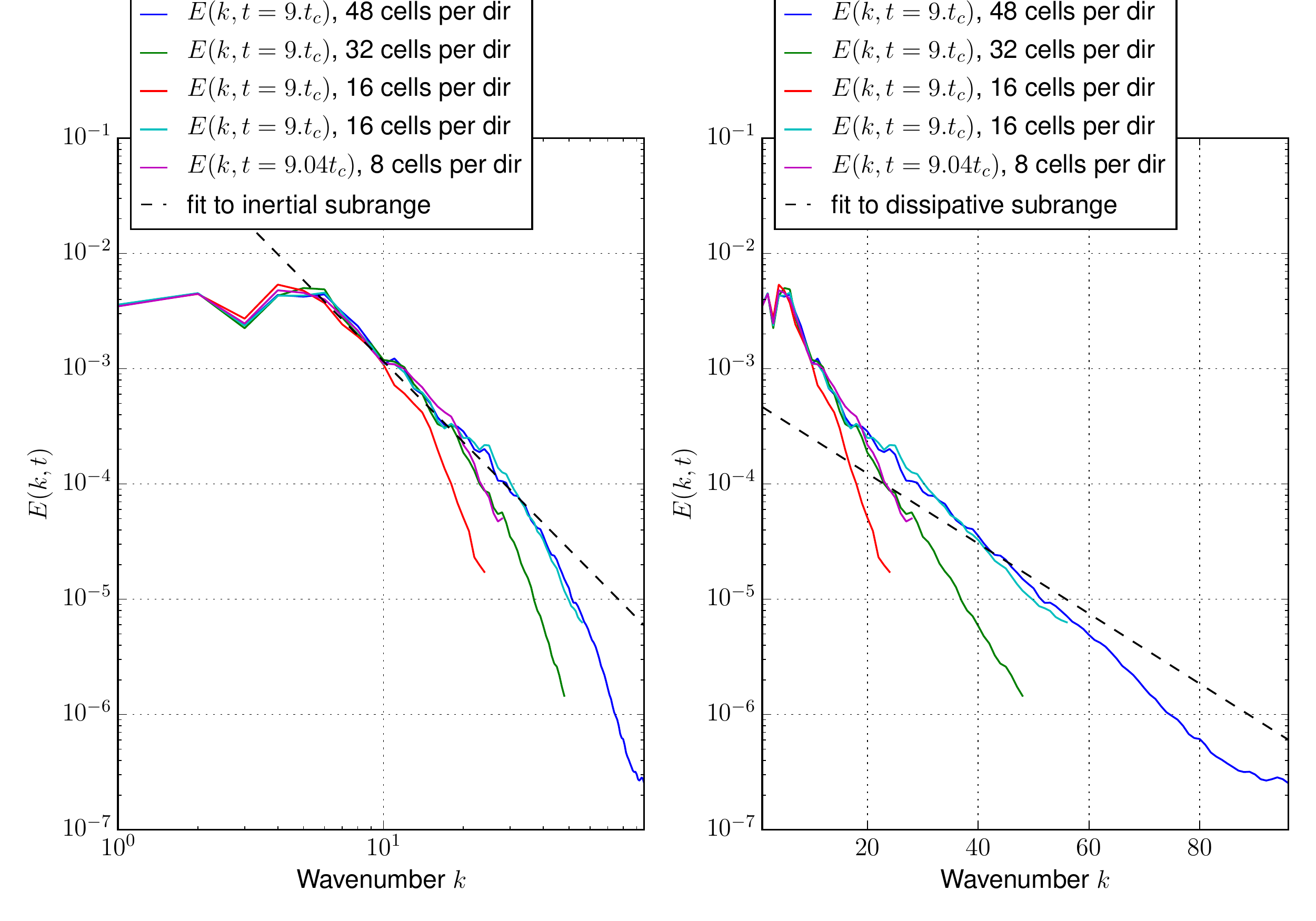}
  \caption{Kinetic spectra at \(t=9t_c\) computed for varying \(p\) and \(h\), but same number of velocity DOFs. The fully
    resolved simulation with \(p=4\) has been added for comparison. All DG configurations here use
    Helmholtz flux reconstruction in \(\RT_h^{p-1}\).}
  \label{fig:taylorgreen3d_kinetic_spectra_samedof}
\end{figure}

The kinetic energy spectra computed at \(t=9t_c\) for the above presented grid configurations are
displayed in the figures \ref{fig:taylorgreen3d_kinetic_spectra_hcoarsen}-~\ref{fig:taylorgreen3d_kinetic_spectra_samedof}.
An approximation quality of a CFD method is to at first reproduce the schematic distribution of
spectral energy caused by the three regimes in a turbulent flow. In explicit, a curve shall exhibit
an integral range at smallest wavenumbers, a power law in wavenumbers representing the inertial
range and a rapid decay at smallest wavenumbers used for the computation.

Figure \ref{fig:taylorgreen3d_kinetic_spectra_hcoarsen} shows the spectral distribution at time
\(t=9t_c\) for polynomial degree \(p=3\) under \(h\)-coarsening. A DNS simulation with 512 velocity
DOFs has been added including fits to the inertial and dissipative subrange, respectively. The
distribution of \(E(k,t)\) should ideally be a decaying power law in the numerically resolved
inertial subrange, smaller scales below the effective resolution should be suppressed as if they
belong to the dissipative subrange. This is indeed the case as largest scales that are resolved by
all configurations, give the same progression at smallest wavenumbers. There are no oscillations for
the underresolved computations in the shorter inertial range. Moreover, the higher modes below each
discretization resolution are suppressed by an exponential drop which can be seen in a log-lin plot
of the same range on the right half. The lower the accuracy is, preferably the more corresponding
wavenumbers get dissipated in the grid convergence results.

Figure \ref{fig:taylorgreen3d_kinetic_spectra_prefine} shows \(E(k,t=9t_c)\) under simultaneous
\(p\)-refinement and \(h\)-coarsening. As in the previous plot, the green curve represents the
spectrum with effective resolution of 256 velocity DOFs for \(p=3\). We observe that the spectra are
almost identical. The \(p=7\) method with about half velocity DOFs closely matches the spectral
distribution of the other shown configurations up to the beginning of the dissipative
subrange. Since the \(p=7\) configuration has approximately half as many velocity DOFs, observe that
fewer wavenumbers are considered for the Fourier transform.

In figure \ref{fig:taylorgreen3d_kinetic_spectra_samedof} the distribution of spectral energy for
configurations with same amount of velocity DOFs is displayed. For comparison, the fully resolved
simulation with \(p=4\) has been added as in the previous plot (both curves in blue color). Similar
to the \(h\)-coarsening study, the underresolved computations preferably have the same kinetic
energy spectrum up to the corresponding numerical dissipation wavenumber, followed by a rapid decay
afterwards. As already observed for the dissipation rate, the result for \(p=7\) with 64 velocity
DOFs closely matches the result for \(p=3\) with twice the velocity DOFs. The configuration \(p=7\)
with 128 velocity DOFs furthermore reproduces the spectral energy of a resolved simulation (\(p=4\),
240 velocity DOFs)

\noindent
Based on the overall results, we conclude that the DG splitting scheme with Helmholtz flux
Raviart-Thomas projection - which evidences to be pressure robust - can be used as a turbulence
model in a large-eddy simulation with no additional modification.

\section{Conclusion and outlook}
\label{sec:conclusion-outlook}

We described a novel postprocessing technique in the projection step of splitting schemes for
incompressible flow that reconstructs the Helmholtz flux in \(H(\text{div})\) and is computed
locally. The resulting velocity field satisfies the discrete continuity equation, and is pointwise
divergence-free. We performed several numerical experiments with the pressure correction scheme
realized by this discrete Helmholtz decomposition, and compared to other discrete Helmholtz
decompositions with respect to accuracy in space. The other variants include a previously introduced
\(H(\text{div})\) postprocessing based on reconstructing the solution to the pressure Poisson
equation, and stabilization-enhanced projections. The \(H(\text{div})\) postprocessing technique
presented in this paper shows to give a high-order DG, pressure robust in space, splitting
scheme. Employing sole reconstruction of the pressure Poisson flux in contrast has been found to
provide a non pressure robust splitting scheme. The results of the numerical tests illustrate that
pressure robust methods outperform non pressure robust methods, especially at high Reynolds
numbers. The main reason is that in the incompressible Euler limit for \(f=0\), the material
derivative is a gradient field which can be handled more appropriately by pressure robust
methods. We have also observed that penalized as div-div projection do not cure a lack of
pressure robustness, but can counteract so in order to give similar results to those obtained by a
pressure robust method. With the 3D Taylor-Green vortex we have included a commonly used benchmark
for underresolved turbulence computations. We compared the stability of the two considered
reconstruction approaches for various resolutions. The results demonstrate that explicit turbulence
modeling in a DG method is not needed, but pressure robustness is a essential ingredient to realize
a robust method for underresolved turbulent incompressible flows.

Apart from the growth of the velocity errors for \(\nu\rightarrow 0\), the numerical results in this article
obtained by div-div projection are equally comparable to Helmholtz flux reconstruction. For the
Helmholtz flux reconstruction in the Raviart-Thomas space of degree \(p-1\), the implementation
currently computes a tentative velocity which is in the DG space of polynomial degree \(p\). A
possible direction for future work is thus to discretize the viscous substep with the anisotropic
tensor product polynomials spanning the local Raviart-Thomas space of degree \(p-1\). However, the
change to anisotropic tensor product polynomials in our spectral DG implementation would produce a
significant amount of additional low-level optimizations leading to sophisticated code which is hard
to maintain. Simultaneous work on a Python-based code generator that can transform a very abstract
description of the variational form into highly optimized code has been shown to give promising
results, \cite{2018arXiv181208075K}.

\section*{Acknowledgements}
\label{sec:acknowledgements}

The authors acknowledge support by the state of Baden-Württemberg through bwHPC.

\bibliographystyle{elsarticle-num}
\bibliography{literature}

\begin{thebibliography}{10}
\expandafter\ifx\csname url\endcsname\relax
  \def\url#1{\texttt{#1}}\fi
\expandafter\ifx\csname urlprefix\endcsname\relax\def\urlprefix{URL }\fi
\expandafter\ifx\csname href\endcsname\relax
  \def\href#1#2{#2} \def\path#1{#1}\fi

\bibitem{KATUL2006117}
G.~Katul, A.~Porporato, D.~Cava, M.~Siqueira,
  \href{http://www.sciencedirect.com/science/article/pii/S0167278906000613}{An
  analysis of intermittency, scaling, and surface renewal in atmospheric
  surface layer turbulence}, Physica D: Nonlinear Phenomena 215~(2) (2006) 117
  -- 126.
\newblock \href {http://dx.doi.org/https://doi.org/10.1016/j.physd.2006.02.004}
  {\path{doi:https://doi.org/10.1016/j.physd.2006.02.004}}.
\newline\urlprefix\url{http://www.sciencedirect.com/science/article/pii/S0167278906000613}

\bibitem{Piatkowski2019}
S.-M. Piatkowski, {A Spectral Discontinuous Galerkin method for incompressible
  flow with Applications to turbulence}, Ph.D. thesis, Heidelberg University
  (2019).
\newblock \href {http://dx.doi.org/https://doi.org/10.11588/heidok.00026674}
  {\path{doi:https://doi.org/10.11588/heidok.00026674}}.

\bibitem{doi:10.1029/162GM07}
A.~H. Al-Yaarubi, C.~C. Pain, C.~A. Grattoni, R.~W. Zimmerman,
  \href{https://agupubs.onlinelibrary.wiley.com/doi/abs/10.1029/162GM07}{Navier-Stokes
  Simulations of Fluid Flow Through a Rock Fracture}, American Geophysical
  Union (AGU), 2013, pp. 55--64.
\newblock \href
  {http://arxiv.org/abs/https://agupubs.onlinelibrary.wiley.com/doi/pdf/10.1029/162GM07}
  {\path{arXiv:https://agupubs.onlinelibrary.wiley.com/doi/pdf/10.1029/162GM07}},
  \href {http://dx.doi.org/10.1029/162GM07} {\path{doi:10.1029/162GM07}}.
\newline\urlprefix\url{https://agupubs.onlinelibrary.wiley.com/doi/abs/10.1029/162GM07}

\bibitem{doi:10.1029/2007GL030545}
M.~B. Cardenas, D.~T. Slottke, R.~A. Ketcham, J.~M. Sharp~Jr., Navier-stokes
  flow and transport simulations using real fractures shows heavy tailing due
  to eddies, Geophysical Research Letters 34~(14).
\newblock \href {http://dx.doi.org/10.1029/2007GL030545}
  {\path{doi:10.1029/2007GL030545}}.

\bibitem{doi:10.1177/1687814019873250}
M.~J. Ahammad, M.~A. Rahman, J.~Alam, S.~Butt, A computational fluid dynamics
  investigation of the flow behavior near a wellbore using three-dimensional
  navier–stokes equations, Advances in Mechanical Engineering 11~(9) (2019)
  1687814019873250.
\newblock \href {http://dx.doi.org/10.1177/1687814019873250}
  {\path{doi:10.1177/1687814019873250}}.

\bibitem{RAEINI20125653}
A.~Q. Raeini, M.~J. Blunt, B.~Bijeljic,
  \href{http://www.sciencedirect.com/science/article/pii/S0021999112001830}{Modelling
  two-phase flow in porous media at the pore scale using the volume-of-fluid
  method}, Journal of Computational Physics 231~(17) (2012) 5653 -- 5668.
\newblock \href {http://dx.doi.org/https://doi.org/10.1016/j.jcp.2012.04.011}
  {\path{doi:https://doi.org/10.1016/j.jcp.2012.04.011}}.
\newline\urlprefix\url{http://www.sciencedirect.com/science/article/pii/S0021999112001830}

\bibitem{doi:10.1002/fld.3653}
F.~Heimann, C.~Engwer, O.~Ippisch, P.~Bastian,
  \href{https://onlinelibrary.wiley.com/doi/abs/10.1002/fld.3653}{An unfitted
  interior penalty discontinuous galerkin method for incompressible
  navier–stokes two-phase flow}, International Journal for Numerical Methods
  in Fluids 71~(3) (2013) 269--293.
\newblock \href
  {http://arxiv.org/abs/https://onlinelibrary.wiley.com/doi/pdf/10.1002/fld.3653}
  {\path{arXiv:https://onlinelibrary.wiley.com/doi/pdf/10.1002/fld.3653}},
  \href {http://dx.doi.org/10.1002/fld.3653} {\path{doi:10.1002/fld.3653}}.
\newline\urlprefix\url{https://onlinelibrary.wiley.com/doi/abs/10.1002/fld.3653}

\bibitem{GAUGER201989}
N.~R. {Gauger}, A.~{Linke}, P.~W. {Schroeder}, {On high-order pressure-robust
  space discretisations, their advantages for incompressible high Reynolds
  number generalised Beltrami flows and beyond}, SMAI Journal of Computational
  Mathematics 5 (2019) 89--129.
\newblock \href {http://dx.doi.org/10.5802/smai-jcm.44}
  {\path{doi:10.5802/smai-jcm.44}}.

\bibitem{LINKE2019350}
A.~Linke, L.~G. Rebholz,
  \href{http://www.sciencedirect.com/science/article/pii/S0021999119301937}{{Pressure-induced
  locking in mixed methods for time-dependent (Navier--)Stokes equations}},
  Journal of Computational Physics 388 (2019) 350--356.
\newblock \href {http://dx.doi.org/10.1016/j.jcp.2019.03.010}
  {\path{doi:10.1016/j.jcp.2019.03.010}}.
\newline\urlprefix\url{http://www.sciencedirect.com/science/article/pii/S0021999119301937}

\bibitem{Schroeder2018}
P.~W. Schroeder, C.~Lehrenfeld, A.~Linke, G.~Lube,
  \href{https://doi.org/10.1007/s40324-018-0157-1}{{Towards computable flows
  and robust estimates for inf-sup stable FEM applied to the time-dependent
  incompressible Navier--Stokes equations}}, SeMA Journal 75~(4) (2018)
  629--653.
\newblock \href {http://dx.doi.org/10.1007/s40324-018-0157-1}
  {\path{doi:10.1007/s40324-018-0157-1}}.
\newline\urlprefix\url{https://doi.org/10.1007/s40324-018-0157-1}

\bibitem{AhmedLinkeMerdon18}
A.~Naveed, L.~Alexander, M.~Christian,
  \href{https://www.degruyter.com/view/j/cmam.2018.18.issue-3/cmam-2017-0047/cmam-2017-0047.xml}{{Towards
  Pressure-Robust Mixed Methods for the Incompressible Navier–Stokes
  Equations}}, Computational Methods in Applied Mathematics 18 (2018) 353--372.
\newblock \href {http://dx.doi.org/10.1515/cmam-2017-0047}
  {\path{doi:10.1515/cmam-2017-0047}}.
\newline\urlprefix\url{https://www.degruyter.com/view/j/cmam.2018.18.issue-3/cmam-2017-0047/cmam-2017-0047.xml}

\bibitem{doi:10.1137/17M1112017}
N.~Ahmed, A.~Linke, C.~Merdon, \href{https://doi.org/10.1137/17M1112017}{{On
  Really Locking-Free Mixed Finite Element Methods for the Transient
  Incompressible Stokes Equations}}, SIAM Journal on Numerical Analysis 56~(1)
  (2018) 185--209.
\newblock \href {http://arxiv.org/abs/https://doi.org/10.1137/17M1112017}
  {\path{arXiv:https://doi.org/10.1137/17M1112017}}, \href
  {http://dx.doi.org/10.1137/17M1112017} {\path{doi:10.1137/17M1112017}}.
\newline\urlprefix\url{https://doi.org/10.1137/17M1112017}

\bibitem{Girault_adiscontinuous}
V.~Girault, B.~Rivière, Mary, F.~Wheeler, A discontinuous galerkin method with
  nonoverlapping domain decomposition for the {S}tokes and {N}avier-{S}tokes
  problems, Math. Comp (2004) 53--84.

\bibitem{girault_riviere}
B.~Rivière, V.~Girault,
  \href{http://www.sciencedirect.com/science/article/pii/S0045782505002690}{Discontinuous
  finite element methods for incompressible flows on subdomains with
  non-matching interfaces}, Computer Methods in Applied Mechanics and
  Engineering 195~(25–28) (2006) 3274 -- 3292, discontinuous Galerkin
  Methods.
\newblock \href {http://dx.doi.org/http://dx.doi.org/10.1016/j.cma.2005.06.014}
  {\path{doi:http://dx.doi.org/10.1016/j.cma.2005.06.014}}.
\newline\urlprefix\url{http://www.sciencedirect.com/science/article/pii/S0045782505002690}

\bibitem{10.2307/4100078}
S.~Zhang, \href{http://www.jstor.org/stable/4100078}{A new family of stable
  mixed finite elements for the 3d stokes equations}, Mathematics of
  Computation 74~(250) (2005) 543--554.
\newline\urlprefix\url{http://www.jstor.org/stable/4100078}

\bibitem{Cockburn2007}
B.~Cockburn, G.~Kanschat, D.~Schötzau,
  \href{https://doi.org/10.1007/s10915-006-9107-7}{{A Note on Discontinuous
  Galerkin Divergence-free Solutions of the Navier–Stokes Equations}},
  Journal of Scientific Computing 31~(1) (2007) 61–73.
\newblock \href {http://dx.doi.org/10.1007/s10915-006-9107-7}
  {\path{doi:10.1007/s10915-006-9107-7}}.
\newline\urlprefix\url{https://doi.org/10.1007/s10915-006-9107-7}

\bibitem{Lehrenfeld2010}
C.~Lehrenfeld, Hybrid discontinuous galerkin methods for solving incompressible
  flow problems, Ph.D. thesis (05 2010).

\bibitem{LINKE2012837}
A.~Linke,
  \href{http://www.sciencedirect.com/science/article/pii/S1631073X1200283X}{A
  divergence-free velocity reconstruction for incompressible flows}, Comptes
  Rendus Mathematique 350~(17) (2012) 837 -- 840.
\newblock \href {http://dx.doi.org/https://doi.org/10.1016/j.crma.2012.10.010}
  {\path{doi:https://doi.org/10.1016/j.crma.2012.10.010}}.
\newline\urlprefix\url{http://www.sciencedirect.com/science/article/pii/S1631073X1200283X}

\bibitem{LINKE2014782}
A.~Linke,
  \href{http://www.sciencedirect.com/science/article/pii/S0045782513002636}{On
  the role of the helmholtz decomposition in mixed methods for incompressible
  flows and a new variational crime}, Computer Methods in Applied Mechanics and
  Engineering 268 (2014) 782 -- 800.
\newblock \href {http://dx.doi.org/https://doi.org/10.1016/j.cma.2013.10.011}
  {\path{doi:https://doi.org/10.1016/j.cma.2013.10.011}}.
\newline\urlprefix\url{http://www.sciencedirect.com/science/article/pii/S0045782513002636}

\bibitem{DIPIETRO2016175}
D.~A.~D. Pietro, A.~Ern, A.~Linke, F.~Schieweck,
  \href{http://www.sciencedirect.com/science/article/pii/S0045782516301189}{A
  discontinuous skeletal method for the viscosity-dependent stokes problem},
  Computer Methods in Applied Mechanics and Engineering 306 (2016) 175 -- 195.
\newblock \href {http://dx.doi.org/https://doi.org/10.1016/j.cma.2016.03.033}
  {\path{doi:https://doi.org/10.1016/j.cma.2016.03.033}}.
\newline\urlprefix\url{http://www.sciencedirect.com/science/article/pii/S0045782516301189}

\bibitem{doi:10.1137/15M1047696}
V.~John, A.~Linke, C.~Merdon, M.~Neilan, L.~G. Rebholz,
  \href{https://doi.org/10.1137/15M1047696}{On the divergence constraint in
  mixed finite element methods for incompressible flows}, SIAM Review 59~(3)
  (2017) 492--544.
\newblock \href {http://arxiv.org/abs/https://doi.org/10.1137/15M1047696}
  {\path{arXiv:https://doi.org/10.1137/15M1047696}}, \href
  {http://dx.doi.org/10.1137/15M1047696} {\path{doi:10.1137/15M1047696}}.
\newline\urlprefix\url{https://doi.org/10.1137/15M1047696}

\bibitem{refId0}
{Linke, Alexander}, {Matthies, Gunar}, {Tobiska, Lutz},
  \href{https://doi.org/10.1051/m2an/2015044}{Robust arbitrary order mixed
  finite element methods for the incompressible stokes equations with pressure
  independent velocity errors}, ESAIM: M2AN 50~(1) (2016) 289--309.
\newblock \href {http://dx.doi.org/10.1051/m2an/2015044}
  {\path{doi:10.1051/m2an/2015044}}.
\newline\urlprefix\url{https://doi.org/10.1051/m2an/2015044}

\bibitem{LINKE2016304}
A.~Linke, C.~Merdon,
  \href{http://www.sciencedirect.com/science/article/pii/S0045782516302730}{Pressure-robustness
  and discrete helmholtz projectors in mixed finite element methods for the
  incompressible navier–stokes equations}, Computer Methods in Applied
  Mechanics and Engineering 311 (2016) 304 -- 326.
\newblock \href {http://dx.doi.org/https://doi.org/10.1016/j.cma.2016.08.018}
  {\path{doi:https://doi.org/10.1016/j.cma.2016.08.018}}.
\newline\urlprefix\url{http://www.sciencedirect.com/science/article/pii/S0045782516302730}

\bibitem{KRANK2017634}
B.~Krank, N.~Fehn, W.~A. Wall, M.~Kronbichler,
  \href{http://www.sciencedirect.com/science/article/pii/S0021999117305478}{A
  high-order semi-explicit discontinuous galerkin solver for 3d incompressible
  flow with application to dns and les of turbulent channel flow}, Journal of
  Computational Physics 348 (2017) 634 -- 659.
\newblock \href {http://dx.doi.org/https://doi.org/10.1016/j.jcp.2017.07.039}
  {\path{doi:https://doi.org/10.1016/j.jcp.2017.07.039}}.
\newline\urlprefix\url{http://www.sciencedirect.com/science/article/pii/S0021999117305478}

\bibitem{FEHN2018667}
N.~Fehn, W.~A. Wall, M.~Kronbichler,
  \href{http://www.sciencedirect.com/science/article/pii/S0021999118304157}{{Robust
  and efficient discontinuous Galerkin methods for under-resolved turbulent
  incompressible flows}}, Journal of Computational Physics 372 (2018) 667--693.
\newblock \href {http://dx.doi.org/10.1016/j.jcp.2018.06.037}
  {\path{doi:10.1016/j.jcp.2018.06.037}}.
\newline\urlprefix\url{http://www.sciencedirect.com/science/article/pii/S0021999118304157}

\bibitem{Schroeder2019}
P.~W. Schroeder, {Robustness of High-Order Divergence-Free Finite Element
  Methods for Incompressible Computational Fluid Dynamics}, Ph.D. thesis,
  University of Göttingen (2019).

\bibitem{doi:10.1002/fld.4763}
N.~Fehn, M.~Kronbichler, C.~Lehrenfeld, G.~Lube, P.~W. Schroeder,
  \href{https://onlinelibrary.wiley.com/doi/abs/10.1002/fld.4763}{High-order dg
  solvers for underresolved turbulent incompressible flows: A comparison of l2
  and h(div) methods}, International Journal for Numerical Methods in Fluids
  0~(0).
\newblock \href {http://dx.doi.org/10.1002/fld.4763}
  {\path{doi:10.1002/fld.4763}}.
\newline\urlprefix\url{https://onlinelibrary.wiley.com/doi/abs/10.1002/fld.4763}

\bibitem{2016arXiv161200657P}
M.~{Piatkowski}, S.~{M{\"u}thing}, P.~{Bastian},
  \href{https://www.sciencedirect.com/science/article/pii/S0021999117308732}{{A
  Stable and High-Order Accurate Discontinuous Galerkin Based Splitting Method
  for the Incompressible Navier-Stokes Equations}}, Journal of Computational
  Physics 356 (2018) 220 -- 239.
\newblock \href {http://dx.doi.org/https://doi.org/10.1016/j.jcp.2017.11.035}
  {\path{doi:https://doi.org/10.1016/j.jcp.2017.11.035}}.
\newline\urlprefix\url{https://www.sciencedirect.com/science/article/pii/S0021999117308732}

\bibitem{FEHN2017392}
N.~Fehn, W.~A. Wall, M.~Kronbichler,
  \href{http://www.sciencedirect.com/science/article/pii/S0021999117306915}{On
  the stability of projection methods for the incompressible navier–stokes
  equations based on high-order discontinuous galerkin discretizations},
  Journal of Computational Physics 351 (2017) 392 -- 421.
\newblock \href {http://dx.doi.org/https://doi.org/10.1016/j.jcp.2017.09.031}
  {\path{doi:https://doi.org/10.1016/j.jcp.2017.09.031}}.
\newline\urlprefix\url{http://www.sciencedirect.com/science/article/pii/S0021999117306915}

\bibitem{doi:10.1002/fld.4683}
N.~Fehn, W.~A. Wall, M.~Kronbichler,
  \href{https://onlinelibrary.wiley.com/doi/abs/10.1002/fld.4683}{A matrix-free
  high-order discontinuous galerkin compressible navier-stokes solver: A
  performance comparison of compressible and incompressible formulations for
  turbulent incompressible flows}, International Journal for Numerical Methods
  in Fluids 89~(3) (2019) 71--102.
\newblock \href
  {http://arxiv.org/abs/https://onlinelibrary.wiley.com/doi/pdf/10.1002/fld.4683}
  {\path{arXiv:https://onlinelibrary.wiley.com/doi/pdf/10.1002/fld.4683}},
  \href {http://dx.doi.org/10.1002/fld.4683} {\path{doi:10.1002/fld.4683}}.
\newline\urlprefix\url{https://onlinelibrary.wiley.com/doi/abs/10.1002/fld.4683}

\bibitem{Braack2014}
M.~Braack, P.~Mucha, Directional do-nothing condition for the navier-stokes
  equations 32 (2014) 507--521.

\bibitem{GiraultRaviartBook}
V.~Girault, P.~Raviart, Finite Element Methods for the Navier-Stokes Equations,
  Springer, 1986.

\bibitem{TemamBook}
R.~Témam, Navier-Stokes Equations. Theory and numerical analysis, North
  Holland, Amsterdam, 1987.

\bibitem{HoustonHartmann2008}
P.~Houston, R.~Hartmann, An optimal order interior penalty discontinuous
  {G}alerkin discretization of the compressible {N}avier-{S}tokes equations, J.
  Comp. Phys. 227 (2008) 9670--9685.

\bibitem{amg4dg}
P.~Bastian, M.~Blatt, R.~Scheichl, Algebraic multigrid for discontinuous
  {G}alerkin discretizations, Numer. Linear Algebra Appl. 19~(2) (2012)
  367--388.
\newblock \href {http://dx.doi.org/10.1002/nla.1816}
  {\path{doi:10.1002/nla.1816}}.

\bibitem{Hillewaert13}
K.~Hillewaert, Development of the discontinuous galerkin method for
  high-resolution, large scale cfd and acoustics in industrial geometries,
  Ph.D. thesis, Univ. de Louvain (2013).

\bibitem{ShuyuWheeler2005}
S.~Sun, M.~F. Wheeler, Symmetric and nonsymmetric discontinuous galerkin
  methods for reactive transport in porous media, SIAM Journal on Numerical
  Analysis 43~(1) (2005) 195--219.
\newblock \href {http://dx.doi.org/10.1137/S003614290241708X}
  {\path{doi:10.1137/S003614290241708X}}.

\bibitem{FLD:FLD373}
L.~J.~P. Timmermans, P.~D. Minev, F.~N. Van De~Vosse, An approximate projection
  scheme for incompressible flow using spectral elements, International Journal
  for Numerical Methods in Fluids 22~(7) (1996) 673--688.

\bibitem{Alexander:1977:DIR}
R.~Alexander, Diagonally implicit {Runge--Kutta} methods for stiff
  {O}.{D}.{E}.'s, SIAM Journal on Numerical Analysis 14~(6) (1977) 1006--1021.

\bibitem{BASTIAN2019417}
P.~Bastian, E.~H. Müller, S.~Müthing, M.~Piatkowski,
  \href{http://www.sciencedirect.com/science/article/pii/S0021999119303973}{Matrix-free
  multigrid block-preconditioners for higher order discontinuous galerkin
  discretisations}, Journal of Computational Physics 394 (2019) 417 -- 439.
\newblock \href {http://dx.doi.org/https://doi.org/10.1016/j.jcp.2019.06.001}
  {\path{doi:https://doi.org/10.1016/j.jcp.2019.06.001}}.
\newline\urlprefix\url{http://www.sciencedirect.com/science/article/pii/S0021999119303973}

\bibitem{Brezzi91}
F.~Brezzi, M.~Fortin, Mixed and Hybrid Finite Element Methods, Springer-Verlag,
  1991.

\bibitem{ErnHdiv2007}
A.~Ern, S.~Nicaise, M.~Vohralík,
  \href{http://www.sciencedirect.com/science/article/pii/S1631073X07004360}{An
  accurate {H}(div) flux reconstruction for discontinuous {G}alerkin
  approximations of elliptic problems}, Comptes Rendus Mathematique 345~(12)
  (2007) 709 -- 712.
\newblock \href
  {http://dx.doi.org/http://dx.doi.org/10.1016/j.crma.2007.10.036}
  {\path{doi:http://dx.doi.org/10.1016/j.crma.2007.10.036}}.
\newline\urlprefix\url{http://www.sciencedirect.com/science/article/pii/S1631073X07004360}

\bibitem{doi:10.1137/16M1089964}
P.~L. Lederer, A.~Linke, C.~Merdon, J.~Schöberl,
  \href{https://doi.org/10.1137/16M1089964}{Divergence-free reconstruction
  operators for pressure-robust stokes discretizations with continuous pressure
  finite elements}, SIAM Journal on Numerical Analysis 55~(3) (2017)
  1291--1314.
\newblock \href {http://arxiv.org/abs/https://doi.org/10.1137/16M1089964}
  {\path{arXiv:https://doi.org/10.1137/16M1089964}}, \href
  {http://dx.doi.org/10.1137/16M1089964} {\path{doi:10.1137/16M1089964}}.
\newline\urlprefix\url{https://doi.org/10.1137/16M1089964}

\bibitem{dune08:1}
P.~Bastian, M.~Blatt, A.~Dedner, C.~Engwer, R.~Kl{\"o}fkorn, M.~Ohlberger,
  O.~Sander, A generic grid interface for parallel and adaptive scientific
  computing. {P}art {I}: {A}bstract framework, Computing 82~(2--3) (2008)
  103--119.
\newblock \href {http://dx.doi.org/10.1007/s00607-008-0003-x}
  {\path{doi:10.1007/s00607-008-0003-x}}.

\bibitem{Bastian2016}
P.~Bastian, C.~Engwer, J.~Fahlke, M.~Geveler, D.~Göddeke, O.~Iliev,
  O.~Ippisch, R.~Milk, J.~Mohring, S.~Müthing, M.~Ohlberger, D.~Ribbrock,
  S.~Turek,
  \href{http://dx.doi.org/10.1007/978-3-319-40528-5_1}{{Hardware-Based
  Efficiency Advances in the EXA-DUNE Project}}, Springer International
  Publishing, Cham, 2016, p. 3–23.
\newblock \href {http://dx.doi.org/10.1007/978-3-319-40528-5_1}
  {\path{doi:10.1007/978-3-319-40528-5_1}}.
\newline\urlprefix\url{http://dx.doi.org/10.1007/978-3-319-40528-5_1}

\bibitem{PAZNER2018344}
W.~Pazner, P.-O. Persson,
  \href{http://www.sciencedirect.com/science/article/pii/S0021999117307830}{{Approximate
  tensor-product preconditioners for very high order discontinuous Galerkin
  methods}}, Journal of Computational Physics 354 (2018) 344–369.
\newblock \href {http://dx.doi.org/10.1016/j.jcp.2017.10.030}
  {\path{doi:10.1016/j.jcp.2017.10.030}}.
\newline\urlprefix\url{http://www.sciencedirect.com/science/article/pii/S0021999117307830}

\bibitem{2018CMAME.341..917A}
M.~{Akbas}, A.~{Linke}, L.~G. {Rebholz}, P.~W. {Schroeder}, {The analogue of
  grad-div stabilization in DG methods for incompressible flows: Limiting
  behavior and extension to tensor-product meshes}, Computer Methods in Applied
  Mechanics and Engineering 341 (2018) 917--938.
\newblock \href {http://arxiv.org/abs/1711.04442} {\path{arXiv:1711.04442}},
  \href {http://dx.doi.org/10.1016/j.cma.2018.07.019}
  {\path{doi:10.1016/j.cma.2018.07.019}}.

\bibitem{FVCA8Benchmark}
\href{https://github.com/FranckBoyer/FVCA8_Benchmark/blob/master/Benchmark.pdf}{Benchmark
  for the fvca8 conference finite volume methods for the stokes and
  navier-stokes equations} (September 2016).
\newline\urlprefix\url{https://github.com/FranckBoyer/FVCA8_Benchmark/blob/master/Benchmark.pdf}

\bibitem{EthierSteinmann1994}
C.~R. Ethier, D.~A. Steinmann, Exact fully {3D} {N}avier--{S}tokes solution for
  benchmarking, Internat. J. Numer. Methods Fluids 19 (1994) 369 -- 375.

\bibitem{10.2307/96892}
G.~I. Taylor, A.~E. Green, \href{http://www.jstor.org/stable/96892}{Mechanism
  of the production of small eddies from large ones}, Proceedings of the Royal
  Society of London. Series A, Mathematical and Physical Sciences 158~(895)
  (1937) 499--521.
\newline\urlprefix\url{http://www.jstor.org/stable/96892}

\bibitem{TaylorGreen3D}
\href{http://www.dlr.de/as/hiocfd}{2nd international workshop on high-order cfd
  methods, cologne, germany} (May 2013).
\newline\urlprefix\url{http://www.dlr.de/as/hiocfd}

\bibitem{fenics:book}
A.~Logg, K.-A. Mardal, G.~N. Wells (Eds.),
  \href{http://dx.doi.org/10.1007/978-3-642-23099-8}{Automated Solution of
  Differential Equations by the Finite Element Method}, Vol.~84 of Lecture
  Notes in Computational Science and Engineering, Springer, 2012.
\newblock \href {http://dx.doi.org/10.1007/978-3-642-23099-8}
  {\path{doi:10.1007/978-3-642-23099-8}}.
\newline\urlprefix\url{http://dx.doi.org/10.1007/978-3-642-23099-8}

\bibitem{2018arXiv181208075K}
D.~{Kempf}, R.~{He{\ss}}, S.~{M{\"u}thing}, P.~{Bastian}, {Automatic Code
  Generation for High-Performance Discontinuous Galerkin Methods on Modern
  Architectures}, arXiv e-prints (2018) arXiv:1812.08075\href
  {http://arxiv.org/abs/1812.08075} {\path{arXiv:1812.08075}}.

\end{thebibliography}

\end{document}